\documentclass[11pt]{article}

\usepackage{amscd}      
\usepackage{amssymb}
\usepackage[intlimits]{amsmath}
\usepackage{xypic}      
\LaTeXdiagrams          
\usepackage[all,v2]{xy}
\xyoption{2cell}
\UseAllTwocells
\xyoption{frame}
\CompileMatrices

\usepackage{theorem}

\newtheorem{prop}{Proposition}[section]
\newtheorem{lem}[prop]{Lemma}

\newtheorem{cor}[prop]{Corollary}
\newtheorem{them}[prop]{Theorem}

\theorembodyfont{\upshape}

\newtheorem{defn}[prop]{Definition}

\newtheorem{rmk}[prop]{Remark}

\newtheorem{numex}[prop]{Example}

\newenvironment{pf}{\begin{trivlist}\item[]{\sc Proof.}}%
            {\nolinebreak $\Box$ \end{trivlist}}
            {\nolinebreak $\Box$ \end{trivlist}}

\newcommand{\noprint}[1]{}

\newcommand{\equal}{=}
\newcommand{\zi}{ \zz}

\newcommand{\B}{{\mathrm{B}}}

\renewcommand{\tilde}{\widetilde}

\newcommand{\toto}{\rightrightarrows}

\newcommand{\upst}{^{\ast}}

\newcommand{\lcom}{_{\scriptscriptstyle\bullet}}
\newcommand{\upcom}{^{\scriptscriptstyle\bullet}}

\newcommand{\XX}{{\mathfrak X}}
\newcommand{\tXX}{\widetilde{{\mathfrak X}}}

\newcommand{\UU}{{\mathfrak U}}

\newcommand{\VV}{{\mathfrak V}}

\newcommand{\zz}{{\mathbb Z}}

\newcommand{\nn}{{\mathbb N}}

\newcommand{\cc}{{\mathbb C}}
\newcommand{\rr}{{\mathbb R}}
\newcommand{\rrr}{{\cal R}}

\newcommand{\aA}{{\cal A}}

\newcommand{\hH}{{\cal H}}

\newcommand{\uU}{{\cal U}}

\newcommand{\del}{\partial}

\newcommand{\id}{\mathop{\rm id}\nolimits}

\newcommand{\smalcirc}{\mbox{\tiny{$\circ $}}}

\newcommand{\ldiag}[1]%
       {\makebox[0cm]{${\scriptstyle#1}\downarrow\phantom{\scriptstyle#1}$}}
\newcommand{\ldiagup}[1]%
       {\makebox[0cm]{${\scriptstyle#1}\uparrow\phantom{\scriptstyle#1}$}}
\newcommand{\rdiag}[1]%
       {\makebox[0cm]{$\phantom{\scriptstyle#1}\downarrow{\scriptstyle#1}$}}
\newcommand{\sediagr}[1]%
       {\makebox[0cm]{$\phantom{\scriptstyle#1}\searrow{\scriptstyle#1}$}}
\newcommand{\nediagr}[1]%
       {\makebox[0cm]{$\phantom{\scriptstyle#1}\nearrow{\scriptstyle#1}$}}
\newcommand{\rdiagup}[1]%
       {\makebox[0cm]{$\phantom{\scriptstyle#1}\uparrow{\scriptstyle#1}$}}
\newcommand{\swdiag}[1]%
       {\makebox[0cm]{$\phantom{\scriptstyle#1}\swarrow{\scriptstyle#1}$}}
\newcommand{\sediag}[1]%
       {\makebox[0cm]{${\scriptstyle#1}\searrow\phantom{\scriptstyle#1}$}}
\newcommand{\nediag}[1]%
       {\makebox[0cm]{${\scriptstyle#1}\nearrow\phantom{\scriptstyle#1}$}}

\newcommand{\iso}{\stackrel{\sim}{\rightarrow}}

\newcommand{\doublearrowstack}[2]%
 {{{{\scriptstyle#1}\atop{\textstyle\longrightarrow}}\atop{{\textstyle\longright
arrow}\atop{\scriptstyle#2}}}}
\newcommand{\rightleftarrowstack}[2]%
 {{{{\scriptstyle#1}\atop{\textstyle\longrightarrow}}\atop{{\textstyle\longlefta
rrow}\atop{\scriptstyle#2}}}}
\newcommand{\leftrightarrowstack}[2]%
 {{{{\scriptstyle#1}\atop{\textstyle\longleftarrow}}\atop{{\textstyle\longrighta
rrow}\atop{\scriptstyle#2}}}}

\newcommand{\overtoparrow}%
{\makebox[0cm]{\beginpicture
\setcoordinatesystem units <.8cm,.4cm> point at 0 0
\setplotarea x from -3 to 3, y from 0 to 1
\setquadratic
\plot -3 0 0 1 3 0 /
\put{\vector(3,-1){0}}[Bl] at 3 0
\endpicture}}

\newcommand{\underbottomarrow}%
{\makebox[0cm]{\beginpicture
\setcoordinatesystem units <.8cm,.4cm> point at 0 0
\setplotarea x from -3 to 3, y from 0 to 1
\setquadratic
\plot -3 1 0 0 3 1 /
\put{\vector(3,1){0}}[Bl] at 3 1
\endpicture}}

\newcommand{\ses}[5]%
{0\longrightarrow#1\stackrel{#2}{ \longrightarrow}#3\stackrel{#4}{
\longrightarrow}#5\longrightarrow0}

\newcommand{\dt}[6]%
{#1\stackrel{#2}{\longrightarrow}#3 \stackrel{#4}{\longrightarrow}#5
\stackrel{#6}{\longrightarrow} #1[1]}

\newcommand{\cat}[1]%
{(\mbox{\rm #1})}

\newcommand{\gm}{\Gamma}
\newcommand{\tgm}{{\tilde{\Gamma}}}
\newcommand{\teta}{{\tilde{\eta}}}

\newcommand{\be }{\begin{eqnarray*}}
\newcommand{\ee }{\end{eqnarray*}}

\def\gpd{\,\lower1pt\hbox{$\longrightarrow$}\hskip-.24in\raise2pt
             \hbox{$\longrightarrow$}\,}

\title{Chern character for twisted $K$-theory of orbifolds}
\author{ Jean-Louis Tu\\
Laboratoire de Math\'ematiques et Applications de Metz\\
Universit\'e de Metz\\
 ISGMP, B\^atiment A, Ile du Saulcy\\
 57000 Metz, France\\
 {\sf email: tu@univ-metz.fr}
 \\\\
Ping Xu \thanks{ Research partially supported by NSF
       grant DMS03-06665 and NSA grant 03G-142. }\\
        Department of Mathematics\\
         Pennsylvania State University \\
         University Park, PA 16802, USA\\
{\sf email: ping@math.psu.edu }}

\begin{document}
\sloppy
\maketitle

\begin{abstract}
For an orbifold $\XX$ and $\alpha \in H^3(\XX, \zz)$, we introduce the twisted
cohomology $H^*_c(\XX, \alpha)$ and prove that the 
non-commutative Chern character of Connes-Karoubi
establishes an isomorphism between the twisted $K$-groups $K_\alpha^* (\XX)
\otimes \cc$ and the 
twisted cohomology $H^*_c(\XX, \alpha)$. This theorem, on the
one hand, generalizes a classical result of Baum-Connes, Brylinski-Nistor, and
others, that if $\XX$ is an orbifold then the Chern character establishes an
isomorphism between the $K$-groups of $\XX$ tensored with $\cc$, and the
compactly-supported cohomology of the inertia orbifold. On the other hand, it
also generalizes a recent result of Adem-Ruan regarding the Chern character
isomorphism of twisted orbifold $K$-theory when the orbifold is a global
quotient by a finite group and the twist is a special torsion class, as well as
Mathai-Stevenson's theorem regarding the Chern character isomorphism of twisted
$K$-theory of a compact manifold.
\end{abstract}

{\small \tableofcontents}

\section{Introduction}

Motivated by mathematical physics and especially string theory, there has been a
great deal of interest in twisted $K$-theory \cite{Witten}.
A mathematically rigorous definition of the $K$-theory of a differentiable stack
twisted by an $S^1$-gerbe was introduced in \cite{TXL04} and some fundamental
properties were also established there. However, just as for the usual
$K$-groups, it is in general a very difficult task to compute the twisted
$K$-groups and very few examples are computed explicitly due to their 
complicated nature.

It is a classical result that for a compact manifold $M$, the Chern character
establishes an isomorphism $K^*(M) \otimes \cc \iso H^*_{dR}(M, \cc)$.
Therefore, in a certain sense, modulo torsion, twisted $K$-groups are isomorphic
to cohomology. It is therefore natural to ask what the twisted cohomology of a
differentiable stack would be, so that the ``Chern character'' would give rise
to an isomorphism. For a torsion class $S^1$-gerbe $\alpha$ over a compact
manifold $M$, little changes and one can show that $K^*_\alpha (M)\otimes \cc$
is isomorphic to $H^*_{dR}(M, \cc)$. However, when the $S^1$-gerbe is of
infinite order, some new phenomena appear, even in the manifold case. The usual
way of defining the $K_0$-group as the Grothendieck group of vector bundles no
longer works. Thus one must use an alternative definition of the ``Chern
character''. Nevertheless, the twisted $K$-groups can be defined as the
$K$-groups of some $C^*$-algebra (see \cite{TXL04}), and therefore one expects
to use techniques of non-commutative differential geometry, and especially
Connes-Karoubi's non-commutative Chern character map \cite{Con85}.

Recall that the non-commutative
Chern character maps the $K$-groups of a smooth subalgebra
(which is a dense topological algebra stable under the holomorphic
functional  calculus) of
a $C^*$-algebra to its periodic cyclic homology. For a compact manifold $M$,
Connes proved \cite{Con85} that the periodic cyclic homology of $C^\infty(M)$ is
isomorphic to the de Rham cohomology of $M$, and therefore the non-commutative
Chern character map indeed generalizes the classical Chern character.

More precisely, for an $S^1$-gerbe $\alpha$ over a differentiable stack $\XX$,
the twisted $K$-groups $K^*_\alpha (\XX)$ are defined to be the $K$-groups of
some $C^*$-algebra $C^*(\XX,\alpha)=C^*(\gm,L)$. It
 contains the subalgebra $C^{\infty}_c (\gm,L)$,
 which is stable under the holomorphic functional calculus when $\gm$
is proper \cite{TXL04}. Here $\gm$ is a Lie groupoid representing $\XX$ and
$L=\tgm\times_{S^1}\cc$ is the associated complex line bundle of the
$S^1$-central extension $\tgm\to \gm$ representing the gerbe $\alpha$. Thus an
essential question is to study the periodic cyclic homology groups
$HP_*(C^{\infty}_c (\gm,L))$. 

In this paper, we confine ourselves to the case when $\XX$ is an orbifold, and
thus $\gm$ an \'etale proper groupoid. The main purpose of the paper is to study
$HP_*(C^{\infty}_c (\gm,L))$. In this case, $S^1$-gerbes are classified by
$H^3(\XX, \zz)$ and therefore $\alpha$ can be considered as an element in
$H^3(\XX, \zz)$.

It is not surprising, as in the classical (i.e., non-twisted) case, that the
inertia orbifold comes into the picture.
This was shown by a classical theorem of Baum-Connes.
In \cite{Bau-Con88}, Baum-Connes proved that when $M$ is a manifold endowed with
a proper action of a discrete group $G$, there is a Chern character isomorphism
$${\mathrm{ch}}: K_i(C^*(M\rtimes G)) \to \oplus_{n\in \nn}
H^{i+2n}_c(\hat{M}/G,\cc)\qquad (i=1,2),$$
where $\hat{M}=\amalg_{g\in G} M^g\times\{g\}$. More generally, it is known that
for any compact orbifold $\XX$ there is a Chern character isomorphism from
$K_i(C^*(\XX))\otimes\cc$ to the $\zz_2$-graded cohomology of the inertia
orbifold $\oplus_{n\in \zz} H^{i+2n}(\Lambda \XX,\cc)$ (\cite{Bry-Nis94, Cra99,
Ade-Rua03} in the global quotient case). See Section~\ref{sec:preliminaries}
below for the definition of inertia orbifolds. In the case when $\XX=M/ G$,
$\Lambda\XX$ is just $\hat{M}/ G$.

In order to study the periodic cyclic homology groups $HP_* (C^{\infty}_c
(\gm,L))$, we introduce the notion of twisted orbifold cohomology $H^*_c (\XX,
\alpha )$. Let $\tXX\to \XX$ be an $S^1$-gerbe over an orbifold $\XX$ with
Dixmier-Douady class $\alpha \in H^3(\XX,\zi)$. Let $\tgm\to \gm\toto M$ be an
$S^1$-central extension representing this gerbe which admits a connection
$\theta$, curving $B$, and curvature $\Omega \in \Omega^3(M)^\gm$. Denote by
$L'\to S\gm$ the induced complex line bundle over the inertia groupoid
$\Lambda\gm \toto S\gm $, which is shown to admit a canonical flat connection.
We denote by $\nabla' : \Omega^*_c( S\gm, L')\to \Omega^{*+1}_c( S\gm, L') $ its
corresponding covariant differential. We define the twisted cohomology groups
(with compact supports) $H^*_c(\XX, \alpha )$ to be the cohomology of the
complex
$$(\Omega^*_c(S\gm, L')^\gm ((u)), \nabla'-2\pi i\Omega u \wedge\cdot),$$
where $u$ is a formal variable of degree $-2$, and $((u))$ are the formal
Laurent series in $u$. When $\XX$ is a smooth manifold, this reduces to the
twisted de Rham cohomology \cite{Mat-Ste04}, which by definition is the
cohomology of the complex $(\Omega^*(M ) ((u)) ,d-u \Omega \wedge\cdot)$. On the
other hand, when $\alpha$ is a torsion class arising from a discrete torsion in
the sense of \cite{Rua03}, this reduces to the twisted orbifold cohomology of
Ruan \cite{Rua03}.

The main result of this paper can be outlined by the following

\begin{them}\label{thm:chern}
Let $\XX$ be an orbifold and $\alpha\in H^3(\XX,\zi)$. Assume that $\tgm \to
\gm\toto M$ is an $S^1$-central extension representing the $S^1$-gerbe
determined by the class $\alpha$, which admits a connection $\theta$, a curving
$B$, and a curvature $\Omega$. Then there are isomorphisms
$$K^*_\alpha(\XX)\otimes\cc\stackrel{{\mathrm{ch}}}{\to} HP_*(C_c^\infty(
\gm,L)) \to H^*_c(\XX,\alpha),$$
where ${\mathrm{ch}}$ denotes the non-commutative Chern character
of  Connes-Karoubi.
\end{them}

Our theorem generalizes, on the one hand, the recent theorem of Mathai-Stevenson
\cite{Mat-Ste04} concerning the non-commutative Chern character
 for twisted $K$-theory of
a compact manifold, and on the other hand, a theorem of Adem-Ruan
\cite{Ade-Rua03} regarding the Chern character of twisted $K$-theory of
orbifolds when the orbifold is a global quotient by a finite group $G$ and the
$S^1$-gerbe is a torsion class induced from a central extension of $G$. 

Note that the periodic cyclic homology of (untwisted) groupoid algebras has been
studied extensively by many authors, including Burghelea \cite{Bur85} in the
case of discrete groups, Feigin-Tsygan \cite{Fei-Tsy87} and
 Nistor \cite{Nis90} in the case of a discrete group acting on a manifold,
Baum-Brylinski-MacPherson \cite{BBM} and   Block-Getzler \cite{BG94}
in the case of a compact Lie group acting on a compact manifold, and
Brylinski-Nistor \cite{Bry-Nis94} and Crainic \cite{Cra99} in the case of
\'etale groupoids.


The paper is organized as follows. Section 2 recalls some basic materials
concerning $S^1$-gerbes over orbifolds including the definition of twisted
$K$-theory of orbifolds. Section 3  introduces twisted cohomology. Section 4
is devoted to the proof of the main theorem by introducing a natural chain map
between the chain complex of the periodic cyclic homology and that of twisted
cohomology.

Note that, besides orbifolds, another important case of differentiable stacks
would be quotient stacks, namely those corresponding to transformation
groupoids. Twisted cohomology and the Chern character for this case will be
discussed in a separate paper.

{\bf Acknowledgments.}
Tu would like to thank Penn State University for the hospitality while part of
the work on this project was being done, and in particular for financial support
from the Shapiro fund.

\section{$S^1$-gerbes over orbifolds}
\label{sec:preliminaries}

In this section we recall a few basic facts concerning orbifolds and
$S^1$-gerbes over orbifolds.

\subsection{Orbifolds}

Roughly speaking, an orbifold is obtained by gluing charts consisting
of manifolds endowed with an action of a finite group. To an orbifold $\XX$,
one can associate an \'etale proper
Lie groupoid $\gm\toto M$, called a presentation of the orbifold.
Note that presentations of an orbifold are
not unique. However, they are uniquely determined
up to Morita equivalence \cite{Moe02, TXL04}.
In other words,  there is a one-to-one correspondence
between Morita equivalence classes of \'etale proper groupoids
and orbifolds. The topological space $|\XX|$ underlying the orbifold $\XX$
is then the orbit space $|\gm|:=M/\gm$.
Recall that a groupoid $\gm\toto M$ is \'etale if both the target and
source maps $t,s:\gm\to M$ are local diffeomorphisms,
and it is proper if the map $(t,s):\gm\to M\times M$ is proper.

By a {\em refinement} of an \'etale proper groupoid $\gm\toto M$, we mean
a pair $(\gm'\toto M',f)$, where  $\gm'\toto M'$ is an \'etale proper
groupoid and $f:\gm'\to \gm$ is an \'etale groupoid morphism which
induces a Morita equivalence.  That is to say,
 (i) $f_0:M'\to M$ is an \'etale map which induces a surjection
$|\gm'|\to |\gm|$, and
(ii) the diagram
$$\xymatrix{
\gm'\dto\rto^{s'\times t'} &M'\times M'\dto\\
\gm\rto^{s\times t} & M\times M}$$
is cartesian.
For instance, if  $(U_i)$  are 
open subsets  of $M$ such that $(U_i/\gm)$ is a cover of $|\Gamma|$,
let $\gm'=\amalg_{i,j} \gm_{U_j}^{U_i}$,
where $\gm_{U_j}^{U_i}=\{g\in\gm|\; s(g)\in U_j\mbox{ and }
t(g)\in U_i\}$.  Then $\gm'$ is endowed with a groupoid structure
with unit space $M'=\amalg_i U_i$, product $(i,j,g)(j,k,h)=(i,k,gh)$
and inverse $(i,j,g)^{-1}=(j,i,g^{-1})$. Then $\gm'\toto M'$,
together with the map $f(i,j,g)=g$, is a refinement of $\gm$.

An \'etale proper  groupoid $\gm\toto M$ is said to be
 \emph{nice} if for all $n$, $\gm_n$ is a disjoint union 
of contractible open subsets.
The following result is proved in \cite[Cor.~1.2.5]{Moe-Pro97}.

\begin{prop}
Any \'etale proper groupoid admits a nice refinement.
\end{prop}

\subsection{$S^1$-gerbes over orbifolds}

An $S^1$-gerbe over an orbifold $\XX$ is represented by
an $S^1$-central extension $S^1\to \tgm\to\gm\toto M$
of a groupoid $\gm\toto M$ representing $\XX$.
More precisely,  if $\gm\toto M$ is an
 \'etale proper Lie groupoid representing an orbifold $\XX$,
there is a one-to-one correspondence between
$S^1$-central extensions of $\gm\toto M$ and
$S^1$-gerbes $\tXX$ over $\XX$ whose restriction
to $M$: $\tXX|_{M}$ admits a trivialization  \cite{BX1, BX2}.
We refer the reader to \cite{BX1, BX2} for 
the general theory between $S^1$-central extensions of
Lie groupoids and $S^1$-gerbes over differentiable stacks.

Below we  recall some basic definitions which
are needed in this paper.

\begin{defn}
Let ${\gm}\toto {M}$ be a Lie groupoid. An {\em $S^1$-central
extension }of $\gm\toto M$ consists of

1) a  Lie groupoid ${\tgm}\toto {M}$, together with a morphism of Lie
groupoids $(\pi,\id):[\tgm\toto M]\to[\Gamma\toto M]$, and

2) a left $S^1$-action on $\tgm$, making $\pi:\tgm\to \Gamma$ a (left)
principal $S^1$-bundle.
\noindent These two structures are compatible in the sense that
$(s\cdot x) (t\cdot y)=st \cdot (xy )$,
for all  $ s,t \in S^{1}$ and $(x, y) \in \tgm\times_M\tgm $.
\end{defn}

The following equivalent definition is quite obvious, which
indicates that this is indeed a generalization of usual
group $S^1$-central extensions.

\begin{prop}
Let $\gm\toto M$ be a Lie groupoid. A Lie groupoid
$\tgm\toto M$ is an $S^1$-central extension of
$\gm\toto M$ if and only if there is a groupoid  morphism
$\pi: \tgm \to \gm$, which is the
identity when being restricted to the unit spaces $M$
  such that  its kernel
$\ker \pi$ is isomorphic to the bundle of groups
$M\times S^1$  and lies  in the center of $\tgm$.
\end{prop}

We now  recall the definition of Morita equivalence of $S^1$-central
 extensions \cite{BX2, TXL04}.

\begin{defn}
We say that  two $S^1$-central extensions
 $S^1\to \tgm\to \gm\toto M$
 and $S^1\to \tgm'\to \gm' \toto M'$ are Morita equivalent if
there exists an  $S^1$-equivariant $\tgm\lcom$-$\tgm'\lcom$-bitorsor
$Z$, by which we mean that $Z$ is a $\tgm\lcom$-$\tgm'\lcom$-bitorsor
endowed with an $S^1$-action such that
$$(\lambda r)\cdot z\cdot r'= r\cdot (\lambda z)\cdot r' = r\cdot
z \cdot
(\lambda r')$$
 whenever $(\lambda, r, r', z)\in S^1\times \tgm\times \tgm'\times Z$ and
the products make sense.
\end{defn}

In the sequel, we will identify an  
 $S^1$-gerbe over an orbifold $\XX$ with the Morita equivalence
class of 
 $S^1$-central extensions of  Lie  groupoids representing  the gerbe.

Denote by ${\mathcal{S}}^1$ (resp. ${\mathcal{R}}$)
the sheaf of $S^1$-valued (resp. $\rr$-valued) smooth functions.
The following is a well-known theorem of Giraud  \cite{giraud}.

\begin{them}[Giraud]
\label{thm:giraud}
Isomorphism classes of $S^1$-gerbes over $\XX$ are in one-to-one
correspondence with $H^2\big(\XX, {\mathcal{S}}^1 \big)$.
\end{them}

The exponential sequence $0\to \zz\to \rrr \to {\mathcal{S}}^1\to 0$
gives rise to a  long exact sequence:

\begin{equation}
\label{eq:exact}
\cdots 
  \to H^2 (\XX, \rrr ) \to H^2(\XX,{\mathcal{S}}^1)
\stackrel{\phi}{\to}
 H^3(\XX,\zz) {\to}   H^3(\XX, \rrr) \to \cdots
\end{equation}

Given an $S^1$-gerbe $\tXX$ over $\XX$, we call the image
 $ \phi ([\tXX])\in H^3(\XX,\zz)$, its {\em Dixmier-Douady}
class, where $[\tXX]\in H^2(\XX,{\mathcal{S}}^1)$
its isomorphism class as in Theorem \ref{thm:giraud}.
Since an orbifold can be represented by an  \'etale proper groupoid,
we have $H^2 (\XX, \rrr )=0$ and $H^3(\XX, \rrr)=0$, and therefore
$H^2(\XX,{\mathcal{S}}^1) \stackrel{\phi}{\to}
 H^3(\XX,\zz)$ is indeed an isomorphism \cite[Proposition~2.22]{TXL04}.
Thus we have

\begin{them}
Isomorphism classes of $S^1$-gerbes over an orbifold  $\XX$
are classified by 
 $H^3\big(\XX, \zz )$.
\end{them}

Thus for an orbifold, Dixmier-Douady classes  completely
classify $S^1$-gerbes over it, which is in general false
for a differentiable stack \cite{TXL04}.
Given an $S^1$-central extension $\tgm\to \gm\toto M$ as above,
  let $L=\tgm\times_{S^1} \cc$ be its associated complex line bundle.
Then $L\to \gm\toto M$
is equipped with   an associative bilinear product
\begin{eqnarray*}
L_g\otimes L_h&\to& L_{gh} \ \ \ \forall (g,h)\in \gm_2\\
(\xi,\eta)&\mapsto&\xi\cdot\eta
\end{eqnarray*}
and an antilinear  involution
\begin{eqnarray*}
L_g&\to& L_{g^{-1}}\\
\xi&\mapsto&\xi^*
\end{eqnarray*}
satisfying the following properties:
\begin{itemize}
\item the restriction of the line bundle to the unit space $M$
is isomorphic to the trivial bundle $M\times\cc\to M$;
\item $\forall \xi,\eta\in L_g$, $\langle \xi,\eta\rangle =
\xi^*\cdot\eta\in L_{s(g)}\cong \cc$ defines a scalar product;
\item $(\xi\cdot\eta)^* = \eta^*\cdot \xi^*$.
\end{itemize}

\begin{numex}
Note that when $\gm$ is a nice \'etale proper groupoid,
then all $S^1$-central extensions are topologically trivial.
Therefore they  are in one-one correspondence
with $S^1$-valued  2-cocycles of the groupoid $\Gamma$, i.e.
differentiable maps $c:\gm_2\to S^1$  satisfying the relation

$$c(g,h)c(gh,k)=c(h,k)c(g,hk), \ \ \ \forall (g, h, k)\in \gm_3, $$
modulo coboundaries, i.e. cocycles of the form $c(g,h)=
b(g)b(h)b(gh)^{-1}$.

The groupoid multiplication on $\tgm$ and the 2-cocycle $c$ are
related by the following equation
\begin{equation}
\label{eq:ex1}
(g, \lambda_1) (h, \lambda_2 )\equal (gh, \lambda_1  \lambda_2 c (g, h)),
\end{equation}
where we identify $\tgm$ with $\gm \times S^1$ by choosing a trivialization.
\end{numex}

\begin{rmk}
Let us explain in the language of groupoids what Ruan
calls discrete torsion in \cite{Rua03}. By definition
\cite[Definition~4.6]{Rua03}, a discrete torsion is a cohomology
class in $H^2(\pi_1^{orb}(\XX),S^1)$. One can
show that the group $\pi_1^{orb}(\XX)$ satisfies the following
universal property: for any discrete group $G$, any generalized
morphism $\gm\to G$ factorizes through the canonical generalized
morphism $\gm \to \pi_1^{orb}(\XX)$. Therefore, any discrete torsion
determines a class in $H^2(\XX,{\mathcal{S}}^1)\cong H^3(\XX,\zz)$,
and a class in $H^3(\XX,\zz)$ comes from a discrete torsion if and only
if it comes from the pull-back of an $S^1$-central extension of
a discrete group $G$ by a generalized morphism $\gm\to G$.
\end{rmk}

\subsection{Inertia groupoid}

Let  $\gm\toto M$ be  a proper and \'etale groupoid representing
an  orbifold $\XX$.
Let $S\gm=\{g\in\gm\vert\; s(g)=t(g)\}$ be the space of
closed loops. Then $S\gm$ is a manifold, on which the natural action of
 $\gm$ by conjugation is smooth. Thus one may form the transformation
groupoid $\Lambda\gm : S\gm\rtimes \gm \toto S\gm$, which  is called
 the {\em inertia groupoid}. Its Morita equivalence class $\Lambda\XX$
is called the {\em inertia orbifold} \cite{Moe02}.

If $S^1\to \tgm\to\gm\toto M$ is an $S^1$-central extension,
then the restriction $\tgm'\to S\gm$
of this  $S^1$-bundle to $S\gm$ is naturally endowed with an action of
$\gm$. To see this,  for any $g\in \gm$, let  $\tilde{g} \in  \tgm$ 
be any of its  lifting.
Then for any  $\gamma\in \tgm'$ such that $s(\gamma)=t(\tilde{g})$,
set
\begin{equation}
\label{eqn:tildeg} 
 \gamma\cdot g=\tilde{g}^{-1}\gamma{\tilde{g}}.
\end{equation}
It is simple to see that this $\gm$-action is well defined, i.e.
independent of the choice of the lifting $\tilde{g}$.
Thus $\tgm'\to S\gm$ naturally  becomes  an $S^1$-bundle
over the inertia orbifold $\Lambda\gm \toto S\gm$
(see also \cite[Lemma 6.4.1]{LU02}). Note that indeed we have
$\tgm'\cong S\tgm$.

\begin{prop}\label{prop:line-bundle}
Let  $\gm\toto M$ be  an \'etale and  proper  groupoid.
Then any  $S^1$-central extension $S^1\to \tgm\to\gm\toto M$
 determines an $S^1$-bundle, and thus a  line bundle $L'$,
over the inertia groupoid $\Lambda\gm \toto S\gm$.
\end{prop}

\begin{numex}
Suppose $\XX=M/G$, where $G$ is a finite group acting on a manifold
$M$ by diffeomorphisms. The inertia groupoid is $\hat M \rtimes G$,
where $\hat M = \cup_{g\in G} M^g\times \{g\}$ and $G$ acts on
$\hat M$ by $(x,g)\cdot \gamma = (x\gamma,\gamma^{-1}g\gamma)$.
 If $S^1\to \tilde{G}\to G$ is an $S^1$-central extension, then
the induced line bundle over
the inertia orbifold
 is a local inner system in the sense of
Ruan (see \cite{Rua03}, \cite[Proposition~4.3.2]{LU02}).
\end{numex}

\subsection{Twisted $K$-theory of orbifolds}
Let $\XX$ be an orbifold and  $\alpha \in H^3(\XX,\zz)$. Choose an $S^1$-central
extension $S^1\to\tgm\to \gm\toto M$ representing
$\alpha$, and denote as above
$L=\tgm\times_{S^1}\cc$. The space $C_c^\infty(\gm,L)$ of
smooth, compactly supported sections of the line bundle $L\to\gm$
is endowed with the convolution product
$$(f_1*f_2)(g)=\sum_{h\in\gm^{t(g)}} f_1(h)\cdot f_2(h^{-1}g),$$
and  the adjoint
$$\xi^*(g)=(\xi(g^{-1}))^*,$$
where $f_1(h)\cdot f_2(h^{-1}g)$ is computed using the product
$L_h\otimes L_{h^{-1}g}\to L_g$.

For all $x\in M$, let $\hH_x$ be the Hilbert space obtained
by completing $C_c^\infty(\gm,L)$ with respect to the scalar
product
$$\langle \xi,\eta\rangle = (\xi^* * \eta)(x)
=\sum_{g\in\gm_x}\langle \xi(g),\eta(g)\rangle .$$
Let $(\pi_x(f))\xi
=f*\xi$ ($f\in C_c^\infty(\gm,L)$, $\xi\in \hH_x$).
Then $f\mapsto \pi_x(f)$ is a $*$-representation of $C_c^\infty(\gm,L)$
in $\hH_x$. The $C^*$-algebra $C^*_r(\gm,L)$ is, by definition, the
completion of $C_c^\infty(\gm,L)$ with respect to the norm:
$\sup_{x\in M}\|\pi_x(f)\|$ \cite{TXL04}. Its Morita equivalence class does not
depend on the choice of the presentation $S^1\to\tgm\to\gm\toto M$,
and therefore its $K$-theory is uniquely determined:

\begin{defn}\cite{TXL04}
The twisted $K$-theory group $K^i_\alpha (\XX)$
is defined as $K_i(C^*_r(\gm,L))$.
\end{defn}

\begin{rmk}
Another way to see the $C^*$-algebra $C^*_r(\gm,L)$
is as follows. The cohomology class $\alpha$ determines
a unique Morita equivalence class of $\gm$-equivariant bundles
of $C^*$-algebras over $M$ satisfying Fell's condition \cite{KMRW98}
such that each fiber is isomorphic to
the algebra of compact operators.
 Denote by $A_\alpha$ one of these $C^*$-algebras, then $C^*_r(\gm,L)$ is Morita equivalent
to the crossed-product algebra $A_\alpha\rtimes_r \gm$.
\end{rmk}

\begin{rmk}
There is another definition of twisted $K$-theory as the Grothendieck
group of twisted vector bundles $K^0_{\alpha,vb}(\XX)$ \cite{Ade-Rua03}.
The group $K^0_{\alpha,vb}(\XX)$ is always zero when $\alpha$ is not
a torsion class. In \cite{TXL04}, it is conjectured that the canonical
map $K^0_{\alpha,vb}(\XX)\to K^0_\alpha(\XX)$ is an isomorphism when
$\alpha$ is torsion. From \cite{TXL04}, this conjecture is known to
be true in some special cases, such as
\begin{itemize}
\item[(a)]
$\XX$ is a compact global
 quotient orbifold $M/G$ (where $G$ is a compact Lie group),
and there exists a twisted vector bundle; or
\item[(b)] $\XX$ is a compact manifold.
\end{itemize}
\end{rmk}

\section{Twisted cohomology}

This section is devoted to the introduction of
twisted cohomology of an orbifold. In  the case
of discrete torsion, this is introduced  by Ruan \cite{Rua03}.
Our definition here is, in a certain sense, a combination
of discrete torsion case with the twisted cohomology
of manifolds.

\subsection{De Rham cohomology}

First, let us recall the definition of the
 de Rham double complex of a
Lie groupoid.
Let $\Gamma\toto M$ be a Lie groupoid.  Define for all
 $p\geq0$
$$\Gamma_{p}=
\underbrace{\Gamma\times_{M}\ldots\times_{M}\Gamma}_{\text{$p$ times}}\,,$$
 {\em i.e.} $\Gamma_{p}$ is the manifold of composable
sequences of $p$ arrows in the groupoid $\Gamma\toto M $ (and $\gm_0
=M$).
We have $p+1$ canonical maps $\gm_{p}\to \gm_{p-1}$ 
giving rise to a              diagram
\begin{equation}\label{sim.ma}
\xymatrix{
\ldots \gm_{2}
\ar[r]\ar@<1ex>[r]\ar@<-1ex>[r] & \gm_{1}\ar@<-.5ex>[r]\ar@<.5ex>[r]
&\gm_{0}\,.}
\end{equation}

Consider  the double complex $\Omega^* (\Gamma\lcom)$:
\begin{equation}
\label{eq:DeRham}
\xymatrix{
\cdots&\cdots&\cdots&\\
\Omega^1(\gm_{0})\ar[u]^d\ar[r]^\partial &
\Omega^1(\gm_{1})\ar[u]^d\ar[r]^\partial
&\Omega^1(\gm_{2})\ar[u]^d\ar[r]^\partial&\cdots\\
\Omega^0(\gm_{0})\ar[u]^d\ar[r]^\partial
&\Omega^0(\gm_{1})\ar[u]^d\ar[r]^\partial
&\Omega^0(\gm_{2})\ar[u]^d\ar[r]^\partial&\cdots
}
\end{equation}
Its boundary maps are $d:
\Omega^{k}( \gm_{p} ) \to \Omega^{k+1}( \gm_{p} )$, the usual exterior
differential of differential forms and $\partial
:\Omega^{k}( \gm_{p} ) \to \Omega^{k}( \gm_{p+1} )$,  the
alternating sum of the pull-back maps of (\ref{sim.ma}).
We denote the total differential by $\delta = (-1)^pd+\del$.
The  cohomology groups of the total complex $C^*_{dR}(\Gamma\lcom)$
$$H_{dR}^k(\Gamma\lcom)= H^k\big(\Omega\upcom(\Gamma\lcom)\big)$$
are called the {\em de~Rham cohomology }groups of $\Gamma\toto M$.

Assume that $\gm\toto M$ is an \'etale proper Lie groupoid
representing an orbifold $\XX$. Since the source and target
maps are local diffeomorphisms, any point $g\in \gm$ induces a 
local diffeomorphism from an open neighborhood of $s(g)$ to
 an open neighborhood of $t(g)$. By
$\Omega^k(M)^\gm$, we denote the space of all differential
$k$-forms on $M$ invariant under such induced actions of
all $g\in \gm$. Note that the space  $\Omega^k(M)^\gm$
is independent of  the presentation $\gm\toto M$, so
we also denote it by  $\Omega^k (\XX)$ interchangeably in
the sequel. Also it is simple to see that $\oplus_k \Omega^k(M)^\gm $
is stable under the  de Rham differential, thus
one obtains a cochain complex  $(\Omega^*(M)^\gm,d)$, whose
cohomology groups are denoted by $H^*_{dR}(M)^\gm$.

\begin{lem}
\label{lem:3.1}
Let $\gm\toto M$ be  an \'etale proper Lie  groupoid.  Consider
the double complex (\ref{eq:DeRham}). Then we have
\begin{enumerate}
\item $E_1^{k, 0}=H^k_{dR}(M)^\gm$;
\item $E_1^{k, p}=0$ if $p\geq 1$.
\end{enumerate}
\end{lem}
\begin{pf}
If $\gm\toto M$ is \'etale, the $q$-th
row  $(\Omega^{q-1}(\gm_p),\del)_{p\ge 0}$
is the complex which computes the differentiable groupoid 
cohomology of $\gm$ with values in the $\gm$-vector bundle 
$\wedge^{q-1}T^*M\to M$.  Thus it is acyclic since
 $\gm\toto M$ is proper \cite{Cra03}.
\end{pf} 


As an immediate  consequence,  we obtain the following

\begin{cor}
Let $\gm\toto M$ be  an \'etale proper Lie  groupoid 
representing  an orbifold $\XX$.
Then 
$$H_{dR}^* (\XX)\cong H^*_{dR}(M)^\gm.$$
\end{cor}

\subsection{Connections and the Dixmier-Douady class}

Let $\tgm\stackrel{\pi}{\to} \gm \toto M$ be an $S^1$-central
 extension representing   an $S^1$-gerbe $\tXX$ over an orbifold $\XX$.
By a {\em pseudo-connection}  \cite{BX1, BX2},
 we mean a pair $(\theta,B)$,
where $\theta\in\Omega^1(\tgm)$ is  a connection  1-form for the bundle
$\tgm\to \gm$, and $B\in\Omega^2(M)$ is  a 2-form.
It is simple to check that $\delta(\theta+B) \in Z^3_{dR}(\tgm\lcom)$
 descends to $Z^3_{dR}(\gm\lcom)$, i.e. there exist unique
$\eta\in\Omega^1(\gm_2)$, $\omega\in\Omega^2(\gm)$ and
$\Omega \in \Omega^3(M)$ such that
$$  \delta(\theta+B)= \pi\upst(\eta+\omega+\Omega).$$
Then $\eta+\omega+\Omega $ is called the
{\em  pseudo-curvature }of the pseudo-connection $\theta+B$.
It is known \cite{BX1, BX2} that 
the class $[\eta+\omega+\Omega]\in H^3_{dR}(\gm\lcom)$ is independent of the
choice of  pseudo-connections.  Under the canonical
homomorphism $H^3(\XX ,\zz)\to H^3(\XX ,\rr)\cong H^3_{dR}(\gm\lcom)$,
 the Dixmier-Douady class of $\tgm$ maps to $[\eta+\omega+\Omega ]$.

Following \cite{Brylinski}, we introduce

\begin{defn}\label{def:connection} \cite{BX2}
Given an $S^1$-central extension $S^1\to \tgm \to \gm \toto M$ of
an   \'etale proper groupoid $ \gm \toto M$,

(i) a connection 1-form $\theta\in \Omega^1(\tgm)$
for the $S^1$-bundle $\tgm\to \gm$, such that
$\del\theta=0$,  is a \emph{connection} on $\tgm\lcom\to \gm\lcom$;

(ii) given $\theta$, a 2-form $B\in\Omega^2(M)$, such that $d\theta=\del B$,
is a {\em curving }on $\tgm\lcom\to \gm\lcom$;

(iii) and given $(\theta,B)$, the 3-form
$\Omega =dB\in  \Omega^3(M )^\gm$ is called
the $3$-{\em curvature }of $(\tgm, \theta ,B)$;

(iv) $\tgm\lcom\to \gm\lcom$ is called  a {\em flat} gerbe, if furthermore
$\Omega=0$.
\end{defn}

It is clear that $\Omega $ is a closed $\gm$-invariant $3$-form, and
thus defines a class in $H^3_{dR}(M)^\gm$.
In fact, by the discussion  above,
it is simple to see that  $[\Omega ]\in H^3_{dR}(M)^\gm$ is the 
image of the Dixmier-Douady
class of the $S^1$-gerbe $\tXX$ over $\XX$   under the canonical
homomorphism 
 $H^3(\XX ,\zz)\to H^3(\XX ,\rr)\cong  H^3_{dR}(\gm\lcom)\cong H^3_{dR}(M)^\gm$.

\begin{rmk}\label{rmk:torsion}
Assume that $\tgm \to \gm\toto M$ is an $S^1$-central extension
whose corresponding Dixmier-Douady class is $\alpha$. Then
the following are equivalent:
\begin{itemize}
\item[(i)] $\alpha$ is a torsion class;
\item[(ii)] there exists a flat connection.
\end{itemize}
\end{rmk}

\begin{lem}\label{lem:theta-phi}
Given an $S^1$-central extension of Lie groupoids
$\tgm\stackrel{\pi}{\to} \gm \toto M$, assume that
 $\theta\in \Omega^1(\tgm)$ is  a connection one-form for the
$S^1$-principal bundle $\pi: \tgm\to \gm$. Suppose we are given a
trivialization $\tgm\cong \gm\times S^1$, and denote by
$\nabla = d+2\pi i \phi$ ($\phi\in \Omega^1(\gm)$) the associated
covariant differential on the line bundle $\tgm\times_{S^1}\cc
\cong \gm\times \cc$. Then
\begin{itemize}
\item[(a)] $\theta = dt + \pi^*\phi$, where $dt$ is the pull-back of the
one-form $dt\in \Omega^1(\rr)$ to $S^1\cong \rr/\zz$.
\item[(b)] Let $c:\gm_2\to S^1$ be the groupoid
 2-cocycle associated to the
central extension $S^1\to \tgm\to \gm$. Then $\del\theta=
\pi^*(\del\phi-\frac{dc}{2\pi i c})$.
\end{itemize}
Therefore $\theta$ is a connection for $\tgm$ if and only if
$\del\phi=\frac{dc}{2\pi i c}$.
\end{lem}

\begin{pf}
(a) is trivial. To show (b), recall that the multiplication
on the groupoid
$\tgm\toto M$ is given by $(g,s)(h,t)=(gh,s+t+\tilde{c}(g,h))
\in \gm\times \rr/\zz$, where $c=e^{2\pi i \tilde{c}}$.
It thus follows that
$$\del\theta|_{((g,s)(h,t))}=\del\phi|_{(g,h)} + ds+dt-
d(s+t+\tilde{c}(g,h))=\del\phi|_{(g,h)}-\frac{dc}{2\pi i c}(g,h).$$
\end{pf}

The following proposition indicates that,
  by passing to a nice refinement,
connections with curving always exist for any $S^1$-central extension over
an  \'etale proper groupoid.

\begin{prop}\label{prop:theta}
Let $\tXX\to \XX$ be an $S^1$-gerbe over an orbifold $\XX$, and 
 $\gm\toto M$  a    nice  \'etale proper groupoid representing
$\XX$. Then $\tXX$ is represented by an $S^1$-central extension
$\tgm\to \gm\toto M$, which admits a connection and curving.
%
\end{prop}
\begin{pf}
Choose a trivialization of the principal bundle $\tgm\to \gm$.
Let $c: \gm_2\to S^1$ be the groupoid 
2-cocycle associated to the
$S^1$-central extension. Consider the equation
$$ \del\phi=\frac{dc}{2\pi i c}=  \frac{1}{2\pi i }d \ln{c}, $$
where $\phi\in \Omega^1 (\gm )$. Since $\del d \ln{c}=d \del  \ln{c}=0$,
the equation above always admits a solution according to \cite{Cra03}.
Thus $\theta = dt + \pi^*\phi$ is a connection for $\tgm\lcom\to 
\gm\lcom$ according
to Lemma \ref{lem:theta-phi}.  Now since $ \del d\phi=d \del \phi=
 \frac{1}{2\pi i }d^2 \ln{c}=0$, again according to \cite{Cra03}
 there exists a $B\in \Omega^2 (M)$ such that $ d\phi = \del B$.
Therefore $d\theta =\pi^* d\phi =\del B$. This shows that
$(\theta , B)$ is  a curving. The curvature is
thus given by $\Omega =dB\in \Omega^3 (M)^\gm$.
\end{pf}

\begin{rmk}\label{rmk:same-coho}
It is simple to see that the forms $-\frac{dc}{c}
\in \Omega^1(\gm_2)$, $-d\phi\in \Omega^2(\gm)$ and
$\Omega \in \Omega^3(M)^\gm $ all represent the same cohomology
class $\alpha  \in H^3(\XX, \rr)$.
\end{rmk}

\begin{prop}
Under the same hypothesis of Proposition  \ref{prop:theta},
 the space of connections is an affine space with underlying vector space
$\Omega^1(M)/\Omega^1(M)^\gm$.
\end{prop}
\begin{pf}
%
>From Lemma~\ref{lem:theta-phi}, it follows that
 $\theta'=dt+\pi^* \phi'$
is another connection if and only if $\del(\phi'-\phi)=0$, i.e.
$\phi'-\phi\in\ker\del$.
Since the complex
$$\Omega^1(M)\stackrel{\del}{\to}
\Omega^1(\gm_1)\stackrel{\del}{\to}
\Omega^1(\gm_2)\stackrel{\del}{\to}\ldots$$
is acyclic, we see that $\ker\del \cong
\Omega^1(M)/\Omega^1(M)^\gm$.
\end{pf}

Assume that $P_0\to X_0$ is an $S^1$-bundle over a  Lie
groupoid $X_1\toto X_0$. Recall that \cite{BX2, LTX}  a {\em connection}
is a connection 1-form $\theta$ on $P_0$ (as an $S^1$-bundle over
$X_0$)  such that $\del\theta=0$. A connection is said to be {\em flat}
if moreover $d\theta =0$. Here  $\del\theta =s^*\theta-t^*\theta \in \Omega^1 (P_1)$, and $P_1\toto P_0$ is  the transformation groupoid
$X_1\times_{ X_0}P_0\toto P_0$. As it was shown in \cite{LTX},
equivalently a connection is a connection 1-form $\theta$ on $P_0$
which is basic with respect to the  action of the pseudo-group of
local  bisections of $X_1\toto X_0$. In particular, if $X_1\toto X_0$ is
an \'etale proper groupoid, the latter is equivalent to that
$\theta $ is invariant under the $X_1$-action. 

We end this subsection with the following

\begin{prop}\label{prop:flat-connection}
Let $\tgm\stackrel{\pi}{\to} \gm \toto M$ be an $S^1$-central extension of
groupoids, which is assumed to admit a connection. Then there exists
 an induced connection on its associated  
$S^1$-bundle $\tgm'\to S\gm$  over the inertia
groupoid  $\Lambda \gm$ constructed as in Proposition \ref{prop:line-bundle}.
It is flat if furthermore the connection admits a curving.

Moreover, if $\gm \toto M$ is an  \'etale proper groupoid,
this induced connection is unique, i.e., is independent of
the choices of connections on $\tgm\stackrel{\pi}{\to} \gm \toto M$.
\end{prop}
\begin{pf}
Assume that $\theta \in \Omega^1 (\tgm )$ is a connection of
the $S^1$-central extension $\tgm\to \gm \toto M$.
We {define} $\theta' \in  \Omega^1 (\tgm ')$ as the pull back  of 
$\theta$ to $\tgm '$.
Note that  $\theta \in \Omega^1 (\tgm )$ is a connection if and only if
it  satisfies the property that for any horizontal path
 $t\mapsto (g_t,h_t)\in\tilde{\gm}_2$, the path $t\mapsto g_th_t$
 is again horizontal.
In order to show that $\theta' \in  \Omega^1 (\tgm ')$  defines
a connection for the $S^1$-bundle  $\tgm'\to S\gm$  over the inertia
groupoid  $\Lambda \gm$,  it suffices to
 show that for any horizontal path $t\mapsto h_t\in
\tilde{\gm}'$ and any path $t\mapsto \gamma_t\in \gm$ such that
$h_t$ and $\gamma_t$ are composable, the path $h_t\cdot \gamma_t$
is again horizontal.
Recall that the action of the element $\gamma_t\in \Gamma$ on
$h_t\in \tilde{\gm}'$ is defined by $h_t\cdot \gamma_t
= \tilde{\gamma}_t^{-1} h_t \tilde{\gamma}_t$, where $\tilde{\gamma}_t$
is a lifting of $\gamma_t$. Let us choose $\tilde{\gamma}_t$ such that
$t\mapsto \tilde{\gamma}_t$ is a horizontal lifting of
$t\mapsto \gamma_t$. It is then clear that
$\mapsto \tilde{\gamma}_t^{-1} h_t \tilde{\gamma}_t$ is a
horizontal path. This  proves the first assertion.

Assume that the connection $\theta \in \Omega^1 (\tgm )$ 
admits a curving. That is, $d\theta =\del B$ for some $B\in \Omega^2 (M)$.
Therefore, we have $d\theta'=\del B_{|S\gm}=0$.
That is, $\tgm'\to S\gm$ is a flat bundle over $\Lambda \gm$.

 Now assume that $\gm \toto M$ is   \'etale and proper. 
Note that any two connections 
 always  differ by $ \pi^* \alpha$, where $\alpha \in \Omega^1 (\gm )$
is a one-form satisfying $ \del \alpha=0$. According  to Lemma 
\ref{lem:3.1}(2), we have $\alpha= \del A$,
where $A\in \Omega^1(M)$.  Since $(\del A)_{|S\gm}=0$,
$\theta'$ is independent of the choice of connections, and
is therefore   canonically defined.
\end{pf}

%
%

\subsection{Twisted cohomology}

We are now ready to introduce the twisted cohomology of an
orbifold.  If $E$ is a vector bundle over an orbifold $\XX$,
 by $\Omega^*(\XX,E)$
 we denote the space of $\gm$-invariant $E$-valued
differential forms  $\Omega^*(M,E)^\gm$, where
$\gm\toto M$ is a presentation of $\XX$. Note that
$\Omega^*(M,E)^\gm$ is independent of the choice of
presentations, so this justifies the notation.
Similarly, $\Omega_c^*(\XX,E)$ denotes the subspace
 of $\Omega^*(M,E)^\gm$ consisting
of differential forms whose supports are  compact in the orbit space
$|\gm|$.

Let $\tXX\to \XX$ be an $S^1$-gerbe over an
 orbifold $\XX$ with Dixmier-Douady class $\alpha \in H^3(\XX,\zi)$.
Let $\tgm\to \gm\toto M$ be an $S^1$-central extension
representing this gerbe, which admits a connection $\theta$, curving $B$
and curvature   $\Omega \in \Omega^3(M)^\gm$.
 Denote by
$L'\to S\gm$  the induced complex line bundle
 over the inertia  groupoid $ \Lambda\gm \toto S\gm $,
  which admits a canonical flat connection according
to Proposition \ref{prop:flat-connection}. By 
 $\nabla' : \Omega^*_c( S\gm, L')\to  \Omega^{*+1}_c( S\gm, L') $
we denote its corresponding covariant differential. It is clear
that $\Omega^*_c( S\gm, L')^\gm$ is stable under $\nabla'$ since
$\nabla'$ is $\gm$-invariant.

\begin{defn}
We define the twisted cohomology groups (with compact supports)
$H^*_c(\XX, \alpha )$
as the cohomology of the complex
$$(\Omega^*_c(\Lambda \XX, L')((u)),
\nabla'-2\pi i\Omega  u \wedge\cdot),$$
where $u$ is a formal variable of degree -2, and
 $((u))$ stands for formal Laurent series in $u$.
\end{defn}

Note that  $\nabla'-2\pi i\Omega  u\wedge\cdot$
is indeed a differential since $(\nabla')^2=0$,
$\Omega\wedge\Omega=0$ and $d\Omega=0$.

\begin{rmk}
\label{rmk:alpha'}
\begin{enumerate}
\item
In the definition above, $2\pi i\Omega$
can be replaced by any non-zero multiple of $\Omega$,
for instance $\Omega$.
Indeed, the complex above is isomorphic to
$(\Omega^*_c(\Lambda\XX,L')((u')),
\nabla'-u'\Omega\wedge\cdot)$ by letting $u= \frac{u'}{2\pi i}$.

\item If $\Omega '\in \Omega^3(M)^\gm $ is another 3-curvature as above,
then $(\Omega^*_c(\Lambda\XX,L')((u)),
\nabla'-2\pi i \Omega u\wedge\cdot)$ and $(\Omega^*_c(\Lambda\XX,L')((u)),
\nabla'-2\pi i \Omega' u\wedge\cdot)$ are isomorphic as chain complexes.
However the isomorphism is not canonical.
The induced map in cohomology is unique up to
an automorphism of the form $\omega\mapsto e^{u\beta}\omega$, where
$\beta\in H^2(\XX, \rr)$. To specify which isomorphism to use,
one needs to specify  the connections. This statement will be
made more precise below (Proposition~\ref{prop:two-connections}).

\item Also, recall that $\Omega^*_c(\Lambda\XX,L')$
is,  by definition,  $\Omega^*_c(S\gm ,L')^\gm$. 
The latter does not depend on the choice of the presentation $\gm\toto M$.
\end{enumerate}

In conclusion, the twisted cohomology groups  defined above
only depend on the orbifold $\XX$ and the class $\alpha \in 
H^3(\XX, \zz)$.  
\end{rmk}

\section{Non-commutative Chern Character}

\subsection{The convolution algebra $\Omega_c^*(\gm,L)$}

Assume that  $\gm\toto M$ is an  \'etale proper groupoid
and $\tgm\to \gm\toto M$ is an $S^1$-central extension.
Let $L\to \gm$ be its associated complex line bundle.
By  $\Omega_c^*(\gm, L)$, we denote the space of
all   compactly supported smooth sections of the
vector bundle $\wedge T^*\gm \otimes L\to \gm$.

\begin{lem}
\label{lem:convolution}
$\Omega_c^*(\gm, L)$ admits a convolution  product, which
makes it into a graded associative algebra.
\end{lem}
\begin{pf}
To see this,  note that
given a vector bundle $E\to \gm$ equipped with an associative
product $E_g\otimes E_h \to E_{gh}$, one can define a convolution
product on the space of compactly supported smooth sections
$C_c^\infty(\gm, E)$ by the formula
$$(f_1 * f_2)(g)=\sum_{h\in \gm^{t(g)}} f_1(h)\cdot f_2(h^{-1}g).$$

Now consider $E=\wedge T^*\gm \otimes L$. Given 
any $(g, h)\in \gm_2$ and
  $\omega_1\otimes\xi_1 \in \wedge T^*_{g}\gm \otimes L_g$,
$\omega_2\otimes\xi_2 \in \wedge T^*_{h}\gm \otimes L_h$,
define
$$(\omega_1\otimes\xi_1 ) \cdot ( \omega_2\otimes\xi_2)=
(r_h \omega_1\wedge l_g \omega_2) \otimes \xi_1 \cdot \xi_2,$$
where $l_g$ and $r_h$ denote the left and right actions
of  $\gm$ on  its exterior product of cotangent bundles,
which are well-defined since $\gm$ is \'etale. It is simple
to see that this indeed defines an associative
product $E_g\otimes E_h \to E_{gh}$.



By  $\Omega_c^*(\gm, L)$, we denote the space of 
the  compactly supported smooth sections of this bundle. Thus it
admits a convolution product, which is denoted by $*$.
\end{pf}

Note that   since $M$ is an open and closed
submanifold  of $\gm$ and the restriction of $L$ to $M$ is the trivial
line bundle, $\Omega^* (M)$ can be naturally considered as
a subspace of $\Omega^*(\gm, L)$. Moreover,
the convolution product between elements in $\Omega^*(M)$
and  elements in $\Omega^*(\gm,L)$ are well defined  using exactly
 the same formula above, even if they may not be  compactly supported.
The following lemma  describes
the convolution product between elements in these two spaces.

\begin{lem}
\label{lem:star}
For any $B\in \Omega^* (M)$ and $\omega \in \Omega^*(\gm, L)$,
we have
\begin{itemize}
\item[(i)] $B* \omega =t^*B\wedge \omega$;
\item[(ii)]  $\omega * B=(-1)^{|B|.|\omega|}s^*B\wedge \omega$;
\item[(iii)] $[B, \omega ]=-\del B\wedge \omega$.
\end{itemize}
\end{lem}
\begin{pf} 
This follows from a straightforward verification.
\end{pf}



Now let $\theta \in \Omega^1 (\tgm)$ be a connection 1-form
on the principal $S^1$-bundle $\tgm\to \gm$, which 
induces a linear connection on  the associated
line bundle $L\to \gm$. By $ \nabla:\Omega^*(\gm,L)\to \Omega^{*+1}(\gm,L)$,
 we denote  its corresponding covariant differential.

\begin{lem}\label{lem:connec-star}
Under the same hypothesis as above, 
$\theta$ is a connection for the $S^1$-central extension
$\tgm\to \gm \toto M$ if and only if the following identities
are  satisfied.
\begin{itemize}
\item[(i)] $\nabla(f_1\ast f_2)=f_1*\nabla f_2 +(\nabla f_1)*f_2, \ \ \forall
f_1, \ f_2\in C_c^\infty(\gm,L)$;
\item[(ii)] $\nabla f = df$ for all $f\in C_c^\infty(M,\cc)$,
where $C_c^\infty(M,\cc)\cong C_c^\infty(M,L)$ is considered
as a subspace  of $C_c^\infty(\gm,L)$ and, similarly,
$\Omega^1_c(M,\cc)$ is a subspace   of $\Omega^1_c(\gm,L)$.
\end{itemize}
\end{lem}
\begin{pf}
We show the direct implication (the converse is proved by working backwards).

Let us show that (i) holds at any point $k_0\in \gm$. Since
$\gm\toto M$ is proper, $F=\gm^{t(k_0)}\cap {\mathrm{supp}}(f_1)$
is finite. Let $\varphi_g\in C_c^\infty(\gm)$ ($g\in F$)
such that the $\varphi_g$'s
have disjoint support and $\varphi_g=1$ in a neighborhood of $g$.
Then the evaluation of both sides of (i) at the point $k_0$ remains
unchanged if we replace $f_1$ by $\sum_g \varphi_gf_1$. By linearity,
we may assume that $f_1$ is of the form
$\varphi_{g_0} f_1$ for some $g_0$, i.e.
 $f_1$ is supported on a small neighborhood of $g_0$. Similarly,
$f_2$ can be assumed to be supported on a small neighborhood
of $h_0$ where $g_0h_0=k_0$.

Let $g\mapsto \tilde{g}$ (resp. $h\mapsto\tilde{h}$)
be a local section of $\tilde{\gm}\to\gm$ around $g_0$ (resp. $h_0$).
By Leibniz' rule, we can assume that $f_1(g)=\tilde{g}$ for $g$ close to $g_0$
and $f_2(h)=\tilde{h}$
for $h$ close to $h_0$.
Let
$t\mapsto(g(t),h(t))\in \gm_2$ be a smooth path with $(g(0),h(0))=
(g_0,h_0)$. Then
\begin{eqnarray*}
(f_1*\nabla f_2)(g_0h_0)((gh)'(0))
&=& 2\pi i \theta(\tilde{h}_0)(f'_2(h_0)(h'(0))) \tilde{g}_0\tilde{h}_0\\
(\nabla f_1*f_2)(g_0h_0)((gh)'(0))
&=& 2\pi i \theta(\tilde{g}_0)(f'_1(g_0)(g'(0))) \tilde{g}_0\tilde{h}_0\\
(\nabla(f_1* f_2))(g_0h_0)((gh)'(0))
&=& 2\pi i \theta(\tilde{g}_0\tilde{h}_0)(m_*(f'_1(g_0)(g'(0)),
f'_2(h_0)(h'(0)))) \tilde{g}_0\tilde{h}_0
\end{eqnarray*}
where $m:\tgm_2\to\tgm$ is the multiplication.
Let us explain, for instance, the notations in the first formula above.
$(f_1*\nabla f_2)(g_0h_0)$ is an element of $T^*_{g_0h_0}\gm\otimes
L_{g_0h_0}$, so its evaluation on the vector $(gh)'(0)$ is an
element of $L_{g_0h_0}=L_{k_0}$.
On the right hand side, the evaluation of $\theta(\tilde{h}_0)\in
T^*_{\tilde{h}_0}\tilde{\gm}$ on the vector
$f'_2(h_0)(h'(0))$ is a scalar; multiplying this scalar by
$\tilde{g}_0\tilde{h}_0\in \pi^{-1}(k_0)\subset\tilde{\gm}\subset L$
yields an element of $L_{k_0}$.

Therefore, we get
\begin{eqnarray*}
\lefteqn{(\nabla(f_1* f_2)-\nabla f_1*f_2-f_1*\nabla f_2)(g_0h_0)((gh)'(0))}\\
&=&-2\pi i (\del\theta)(\tilde{g}_0,\tilde{h}_0)
(f'_1(g_0)(g'(0)),f'_2(h_0)(h'(0)))=0.
\end{eqnarray*}

To show (ii), note that (i) obviously holds whenever $f_j\in C^\infty(M)$
(not necessarily compactly supported). For $f_1=f_2=1_M$, we have 
$\nabla(1_M)=0$. Then, by Leibniz' identity, we get
$\nabla(f)=\nabla(f.1_M)=df 1_M + f \nabla(1_M)=df$.
\end{pf}

As a consequence, we have

\begin{prop}
$\nabla:\Omega^*_c(\gm,L)\to \Omega^{*+1}_c(\gm,L)$ satisfies the
Leibniz identity

\begin{equation}
\label{eq:Leibniz}
\nabla(\omega_1*\omega_2)=(\nabla\omega_1)*\omega_2+
(-1)^{|\omega_1|}\omega_1*(\nabla\omega_2).
\end{equation}
\end{prop}
\begin{pf}
It is simple to see,  by using Lemma  (\ref{lem:star}),
that Eq. (\ref{eq:Leibniz}) holds  when
 $\omega_1\in C^\infty_c(M,\cc)$ and  $\omega_2\in \Omega_c^*(\gm,L)$.

To prove Eq. (\ref{eq:Leibniz}) in general, it suffices to  prove
it locally. By linearity, we may assume that for 
$j=1, 2$,  $\omega_j=\eta_j *f_j$
($\eta_j\in\Omega_c^*(M,\cc)$, $f_j\in C_c^\infty(\gm,L)$),
with $f_j$ supported on a small neighborhood of $g_j$ and
$\eta_j$ supported on a small neighborhood of $t(g_j)$, where
$(g_1,g_2)\in \gm_2$. Since $\gm$ is \'etale, there is
a canonical local bisection ${\cal L}$ of $\gm$ through  the point $g_1$.
 It is clear that ${\cal L}$ induces  a local diffeomorphism on $M$, denoted
by $Ad_{{\cal L}^{-1}}$,  from a
small neighborhood of $t(g_1)$ to   a small neighborhood of $s(g_1)$.
>From Lemma \ref{lem:star}, it follows that 
\begin{equation}
\label{eqn:anticom}
\omega*\eta = (-1)^{|\omega|| \eta|} \teta*\omega,
\end{equation}
for any  $\eta \in  \Omega_c^*( M, \cc )$ supported around $s(g_1)$,
 and $\omega\in\Omega_c^*(\gm,L)$ supported around $g_1$, where
$\teta =Ad^*_{{\cal L}^{-1}} \eta$.  Therefore, we have
\begin{eqnarray*}
\nabla(\omega_1*\omega_2) &=&
\nabla(\eta_1*\teta_2*f_1*f_2)\\
&=&d(\eta_1*\teta_2)*f_1*f_2 + (-1)^{|\omega_1|+|\omega_2|}
\eta_1*\teta_2*\nabla(f_1*f_2)\\
&=&(d\eta_1)* \teta_2*f_1*f_2 + (-1)^{|\omega_1|}\eta_1
*d\teta_2* f_1*f_2\\
&&+(-1)^{|\omega_1|+|\omega_2|}
\eta_1* \teta_2*(\nabla f_1)*f_2 +(-1)^{|\omega_1|+|\omega_2|}
\eta_1 *\teta_2* f_1*\nabla f_2\\
&=&(d\eta_1)* f_1*\eta_2 *f_2 + (-1)^{|\omega_1|}\eta_1*
\tilde{d\eta_2}* f_1*f_2\\
&&+(-1)^{|\omega_2|}
\eta_1 *(\nabla f_1)*\eta_2*f_2+(-1)^{|\omega_1|+|\omega_2|}
\eta_1*f_1 *\eta_2*\nabla f_2\\
&=& (\nabla\omega_1)*\omega_2 + (-1)^{| \omega_1|}\omega_1
*\nabla\omega_2.
\end{eqnarray*}
Eq. (\ref{eq:Leibniz})  thus follows.
\end{pf}

\subsection{The trace map $Tr:\Omega^*_c(\gm,L)\to\Omega_c^*(\Lambda\XX,L')$}

This subsection is devoted to the introduction of
the trace map.  Let $\omega\in\Omega^*_c(\gm,L)$.
 For any  $h\in S\gm$ and $g\in \gm$ such
that $s(h)=t(g)$,   $\omega_{g^{-1}hg}\in 
\wedge T^*_{g^{-1}hg}  \gm \otimes L_{g^{-1}hg}$.
Then $i^* \omega_{g^{-1}hg}\in \wedge T^*_{g^{-1}hg} S\gm \otimes L_{g^{-1}hg}
\cong \wedge T^*_{g^{-1}hg} S\gm \otimes  L'_{g^{-1}hg}$, where
$i: S\gm \to \gm$ is the inclusion. Therefore
$g_{*} i^* \omega_{g^{-1}hg}\in  \wedge T^*_h S\gm \otimes  L'_h$. Here
$\wedge T^*  S\gm \otimes L'\to S\gm$ is considered as a vector bundle
over the inertia groupoid $\Lambda \gm\toto S\gm$ and therefore
admits a $\gm$-action, which is denoted by $g_*$. Now  define 
$Tr(\omega)\in \Omega_c^*(S\gm , L')$ by
$$Tr(\omega)_h :=\sum_{g\in\gm^{s(h)}} g_*  i^* \omega_{g^{-1}hg}, \ \ \forall
h\in S\gm  $$
It  is simple to see that $Tr(\omega)$ is  a $\gm$-invariant element
in $\Omega^*(S\gm,L')$ whose support is compact in
$|\Lambda\gm|$. Therefore it defines
   an element in $\Omega_c^*(\Lambda\XX,L')$.

The proposition below describes an important property of this
trace map.

\begin{prop}\label{prop:trace}
\begin{itemize}
\item[(a)] $\forall \omega_1\in \Omega^{|\omega_1|}(\gm, L) $ and
$ \omega_2\in \Omega^{|\omega_2|}(\gm, L)$, we have
$$Tr(\omega_1*\omega_2)=(-1)^{|\omega_1||\omega_2|} Tr(\omega_2*\omega_1).$$
\item[(b)] Assume that $\nabla$ is a linear connection on $L\to \gm$
induced from a connection on the $S^1$-central extension
$\tgm\to\gm \toto M$. Let $\nabla'$ be its induced
 connection   on the line bundle $L'\to S\gm$ over
 $\Lambda \gm \toto S\gm $. Then
$$\nabla'(Tr\omega) = Tr(\nabla\omega). $$
 \end{itemize}
\end{prop}
\begin{pf}
It is simple to see, from a straightforward computation, that
for any $h\in S\gm $,
\begin{eqnarray*}
Tr(\omega_1*\omega_2)_h&=&
\sum_{g\in \gm^{s(h)}, k_1k_2=g^{-1}hg}
g_* i^* (\omega_1(k_1) \cdot \omega_2(k_2))\\
Tr(\omega_2*\omega_1)_h&=&
\sum_{g'\in \gm^{s(h)}, k_2k_1={g'}^{-1}hg'}
g'_* i^* (\omega_2(k_2) \cdot  \omega_1 (k_1)) \ \ \ (\mbox{let } g=g'k_2)\\
&=&\sum_{g\in \gm^{s(h)}, k_1k_2=g^{-1}hg}
(gk_2^{-1})_*i^* (\omega_2(k_2) \cdot  \omega_1(k_1)).
\end{eqnarray*}
It thus suffices to prove the following identity:

\begin{equation}
\label{eq:12}
i^* (\omega_1(k_1)\cdot \omega_2(k_2))=
(-1)^{|\omega_1||\omega_2|} i^* (k_2^{-1})_* (\omega_2(k_2) \cdot \omega_1(k_1)).
\end{equation}

By linearity, we can assume that $\omega_j=\eta_j\otimes\xi_j
\in \wedge T^*_{k_j}\gm \otimes L_{k_j}, \ j=1, 2$. Thus Eq. (\ref{eq:12})
reduces to
\begin{eqnarray}
\label{eqn:com1}
\xi_1\cdot \xi_2 &=& k_{2*}^{-1} (\xi_2\cdot\xi_1) \ \ \mbox{and}\\
\label{eqn:com2}
i^* (\eta_1 \cdot \eta_2) &=&(-1)^{|\eta_1||\eta_2|}
i^* k_{2*}^{-1}(\eta_2\cdot  \eta_1).
\end{eqnarray}
Since $\xi_2\in  L_{k_2}$,
according to Eq. (\ref{eqn:tildeg}),  we have
$$k_{2*}^{-1} (\xi_2\cdot\xi_1) =\xi_2^{-1} \cdot(\xi_2\cdot\xi_1)\cdot
 \xi_2=\xi_1\cdot \xi_2.$$
Thus Eq.  (\ref{eqn:com1}) follows.

For Eq. (\ref{eqn:com2}), note that
\begin{eqnarray*}
&&k_{2*}^{-1}(\eta_2\cdot  \eta_1)\\
&=&k_{2*}^{-1}(r_{k_1}\eta_2 (k_2)\wedge l_{k_2}\eta_1 (k_1) )\\
&=&Ad_{k_2^{-1}} r_{k_1}\eta_2 (k_2)\wedge Ad_{k_2^{-1}}l_{k_2}\eta_1 (k_1) \\
&=&(-1)^{|\eta_1||\eta_2|} (r_{k_2} \eta_1 (k_1) \wedge Ad_{(k_1 k_2)^{-1} }
l_{k_1}  \eta_2 (k_2)\\
&=&(-1)^{|\eta_1||\eta_2|} (r_{k_2} \eta_1 (k_1) \wedge (k_1 k_2)_*^{-1}
(l_{k_1}  \eta_2 (k_2)).
\end{eqnarray*}
Therefore it suffices to prove that 
$$ (k_1 k_2)_*^{-1} i^* l_{k_1}  \eta_2 (k_2) =i^* l_{k_1}  \eta_2 (k_2),$$
since $(k_1 k_2)_*^{-1}$ and $i^*$ commute.
The latter holds due to the following general fact: for any $\gamma \in S\gm$,
and $\eta \in \wedge T^*_\gamma(S\gm)$, one has
\begin{equation}
\label{eqn:heta}
\gamma_*\eta = \eta .
\end{equation}
 Locally, we can assume that $\gm\toto M$
is the crossed product $M\rtimes H\toto M$ of a manifold
by a finite group. We then have $\eta\in \wedge T^*_x M^h$ and  $r=(x, h)$,
where $x\in M^h $ and $h\in H$. Thus Eq. (\ref{eqn:heta}) becomes
$$h_*\eta = \eta ,$$
which is obvious since $h$ acts  on $T_x M^h$  as the identity.
\end{pf}

\subsection{Chain map}

Now we are ready to prove the main theorem.
First we will introduce a chain map
from the chain complex of periodic cyclic cohomology
to the  chain complex of the  twisted cohomology.

>From now on, we assume that   $\tgm\to \gm\toto M$ is an 
$S^1$-central extension of a proper \'etale groupoid
$\gm\toto M$, which admits a connection $\theta$ (with
the corresponding covariant differential $\nabla$ on
the associated line bundle $L\to \gm$), curving $B$
and curvature $\Omega$ as in Definition \ref{def:connection}.
 We start by constructing a chain map
from the cyclic bi-complex to the complex
$(\Omega^*(\Lambda\XX,L'))((u))$.

Let $\aA=C_c^\infty(\gm,L)$, and
then  $\aA^{\otimes n} = C_c^\infty(\gm^n,L^{\otimes n})$ (here we 
use the inductive tensor product, for which we have  $C_c^\infty(M)\otimes
C_c^\infty(N)\cong C_c^\infty(M\times N)$). Denote by $\tilde\aA$
the unitization of $\aA$.

Let $CC_k(\aA)=\tilde{\aA}\otimes \aA^{\otimes k}$.
Recall \cite{BG94} \cite[Prop. 2.2.16]{Lod98}
 that the periodic cyclic homology of $\aA$ is the
cohomology of the complex $(CC_k(\aA)((u)),b+u\B)$,
where $u$ is a formal variable of degree -2, and
$b$ is a differential of degree $-1$ while $\B$ is a 
differential of degree $1$ defined by
\begin{eqnarray*}
b(\tilde{a}_0\otimes a_1\otimes\cdots \otimes  a_k)
&=&\sum_{j=0}^{k-1} (-1)^j \tilde{a}_0\otimes\cdots\otimes a_ja_{j+1}\otimes
\cdots\otimes a_k\\
&&\qquad + (-1)^k a_k\tilde{a}_0\otimes a_1\otimes\cdots\otimes a_{k-1}, \\
\B(\tilde{a}_0\otimes a_1\otimes\cdots a_k)&=&
\sum_{i=0}^k (-1)^{ik}1\otimes a_i\otimes\cdots\otimes a_k
\otimes a_0\otimes\cdots\otimes a_{i-1}.
\end{eqnarray*}

Then $b$ and $\B$ satisfy the identity $b^2=\B^2=b\B+\B b=0$.

Following \cite{Gor99} (see also  \cite{Mat-Ste04}), we  introduce
a linear map
$\tau_{\nabla,B}: CC_k(\aA)\to \Omega^*_c(\Lambda \XX,L')((u))$ by

\begin{equation}
\label{eq:tau}
\tau_{\nabla,B}(\tilde{a}_0\otimes a_1\otimes \cdots\otimes a_k)
=\int_{\Delta^k} Tr(\tilde{a}_0*
e^{ 2\pi i u s_0  B }*\nabla(a_1)*\cdots *
\nabla(a_k)*e^{2\pi i u s_k  B })ds_1\cdots ds_k,
\end{equation}
where the integration is over the $k$-simplex
$\Delta^k =\{(s_0, \cdots, s_k)|s_0\geq 0, \ldots, s_k\geq 0, 
s_0+ \cdots +s_k =1\}$. 
Here  the curving  2-form $B\in \Omega^2(M)$ is considered
as an element of $\Omega^2(\gm,L)$ since $M$ is an open and closed
submanifold  of $\gm$ and the restriction of $L$ to $M$ is the trivial
line bundle.
We then extend $\tau_{\nabla,B}$ to a $\cc((u))$-linear map
$$\tau_{\nabla,B}: \ CC_k(\aA)((u)) \to \Omega^*_c(\Lambda \XX,L')((u)). $$

\begin{prop}
\label{prop:chain-map}
$\tau_{\nabla,B}\smalcirc(b+u \B) = (u\nabla' - 2\pi i u^2 \Omega )\smalcirc\tau_{\nabla,B}$.
\end{prop}
\begin{pf}
First, note that, according to Lemma \ref{lem:star},  we have 
\begin{equation}
\label{eqn:Ba}
\nabla^2 a= 2\pi i \del B \wedge a= -2\pi i [B,a], \ \ \ \forall a\in\Omega^*(\gm,L),
\end{equation}
where $B\in \Omega^2(M,\cc)$ is considered as an element of $\Omega^2(\gm,L)$.



The rest of the  proof is similar to \cite[Proposition~5.4]{Mat-Ste04}.
We will  sketch   the  main steps below.
 Using Proposition~\ref{prop:trace}(b),
the relation $dB=\Omega$ and Eq. (\ref{eqn:Ba}),
 we find that
$(u\nabla '- 2\pi i u^2 \Omega )\smalcirc\tau_{\nabla,B}
(\tilde{a}_0\otimes a_1\otimes \cdots\otimes a_k)$ is the sum of
the following  two terms:
\begin{equation}
\label{eqn:1st-term}
u\int_{\Delta^k} Tr(\nabla(a_0) *e^{2\pi i u s_0 B}*\nabla(a_1)*
\cdots*\nabla(a_k)*e^{2\pi i u s_k B})ds_1\cdots ds_k,
\end{equation}
and
\begin{eqnarray}
\label{eqn:2nd-term}
&&\sum_{i=1}^k (-1)^{i-1}u\int_{\Delta^k} Tr(\tilde{a}_0 *
e^{2\pi i u s_0 B}*\nabla(a_1)* \cdots  *e^{2\pi i u s_{i-1}B}*
[-2\pi i B,a_i]*e^{2\pi i u s_i B}\\
\nonumber
&&\qquad
*\nabla(a_{i+1})*\cdots *\nabla(a_k)*e^{2\pi i u s_k B})
ds_1\cdots ds_k
\end{eqnarray}

Using Proposition~\ref{prop:trace}(a), one shows that the first term
as given by Eq.  (\ref{eqn:1st-term}) is equal to $(\tau_{\nabla,B}\smalcirc u \B)
(\tilde{a}_0\otimes a_1\otimes \cdots\otimes a_k)$.

Using the formula $[a_i,e^K]=\int_0^1 e^{(1-s)K}[a_i,K]e^{sK}\,ds$
 (see \cite{Quillen}), one identifies the second term
  as given by Eq. (\ref{eqn:2nd-term})  with
$(\tau_{\nabla,B}\smalcirc b)(\tilde{a}_0\otimes a_1\otimes \cdots\otimes
a_k)$.
\end{pf}

As a consequence, the chain map $\tau_{\nabla, B}$ induces a homomorphism
 in cohomology
$$\tau_{\nabla,B}: HP_*(\aA)\to H_c^*(\XX,\alpha).$$

In the  next proposition, we study   how $\tau_{\nabla,B}$
depends on the choice of $\nabla$ and $B$.
Assume that $\tgm\stackrel{\pi}{\to} \gm\toto M$   admits 
 connections $\theta_i$,
 curvings $B_i$
and  $3$-curvatures $\Omega_i$, $i=1, 2, 3$, as in
Definition~\ref{def:connection}. Their corresponding
covariant differentials are denoted by  $\nabla_i$, $i=1, 2, 3$.
Then, according to Lemma \ref{lem:3.1},
 $ \theta_1-\theta_0=\del A$ for some $A\in  \Omega^1 (M)$.
Thus,
$$\nabla_1-\nabla_0 =2\pi i\del A .$$
Let 
$$\beta_{(\nabla_0,B_0),(\nabla_1,B_1)} =B_1-B_0-dA. $$ 

\begin{lem}
\label{lem:beta}
\begin{itemize}
\item[(a)] $d\beta_{(\nabla_0,B_0),(\nabla_1,B_1)}=\Omega_1-\Omega_0$;
\item[(b)] 
the form $\beta_{(\nabla_0,B_0),(\nabla_1,B_1)}
\in \Omega^2(M)$ is $\gm$-invariant, and its
class in $\Omega^2(\XX)/d\Omega^1(\XX)$ is well defined; and
\item[(c)] $\beta_{(\nabla_0,B_0),(\nabla_2,B_2)}
=\beta_{(\nabla_0,B_0),(\nabla_1,B_1)}
+\beta_{(\nabla_1,B_1),(\nabla_2,B_2)}$ modulo
$d\Omega^1(\XX)$.
\end{itemize}
\end{lem}
\begin{pf}
(a) and (c) are obvious. 

(b)  Let $\beta =\beta_{(\nabla_0,B_0),(\nabla_1,B_1)}$.
Then  $\del\beta = \del(B_1-B_0)-d\del A =
 \del(B_1-B_0)-d(\theta_1-\theta_0)
=0$. Hence $\beta$ is $\gm$-invariant. Moreover,
if $A'\in \Omega^1 (M) $ is another one-form such that
$ \theta_1-\theta_0=\del A'$,
then $\beta'-\beta = d(A-A')\in d\Omega^1(M)^\gm$ since $\del (A-A')=0$.
\end{pf}

\begin{prop}\label{prop:two-connections}
Let $u'=2\pi i u$. Consider the diagram
\begin{equation}
\label{eq:12a}
\xymatrix{
HP_*(\aA)\ar[r]^{\hspace{-40pt}\tau_{\nabla_0,B_0}}
         \ar[dr]_{\tau_{\nabla_1,B_1}}
  & \qquad H^*(\Omega^*_c(\Lambda\XX,L')((u')),\nabla'-u'\Omega_0)
     \ar[d]^{e^{u'\beta_{(\nabla_0,B_0),(\nabla_1,B_1)}}}\\
  &\qquad H^*(\Omega^*_c(\Lambda\XX,L')((u')),\nabla'-u'\Omega_1)
}
\end{equation}
Then
\begin{itemize}
\item[(a)] the vertical map is well defined, i.e. it is
independent of  the choice of $\beta_{(\nabla_0,\Omega_0),(\nabla_1,\Omega_1)}$
modulo $d\Omega^1(\XX)$; and
\item[(b)] the diagram (\ref{eq:12a}) is commutative.
\end{itemize}
\end{prop}
\begin{pf}
Let $\beta =\beta_{(\nabla_0, B_0),(\nabla_1, B_1)}$.
First,  note that $e^{u'\beta}$ is indeed a chain map, since
$\forall \omega \in \Omega^*_c (S\gm, L')$, we have 
\begin{eqnarray*}
(\nabla'- u'\Omega_1)(e^{u'\beta}\omega)
  &=& e^{u'\beta}(\nabla'\omega + u'
d\beta \omega -u'\Omega_1\omega)\\
  &=& e^{u'\beta}(\nabla' - u'\Omega_0)\omega,
\end{eqnarray*}
according to Lemma \ref{lem:beta} (a).

(a) It suffices to show 
that for all $A\in \Omega^1(M)^\gm$
and $\omega\in \Omega^*_c(\Lambda\XX,L')((u'))$,
$e^{u' dA}\omega$ and $\omega$ are  cohomologous in
$H^*(\Omega^*_c(\Lambda\XX,L')((u')),\nabla'-u'\Omega_0)$.

Let $f(x)=\frac{e^x-1}{x}$, and define 
$$K:\Omega^*_c(\Lambda\XX,L')((u'))
\to \Omega^{*-1}(\Lambda\XX,L')((u'))$$ by $K(\omega)=u'A f(u'dA)\omega$.
Then

\begin{eqnarray*}
(\nabla'-u'\Omega_0 )K(\omega) &=&
u'dA f(u'dA)\omega -u'A f(u'dA)\nabla'\omega -u'^2 \Omega_0 A f(u'dA)\omega\\
&=& (e^{u'dA}-1)\omega -u'Af(u'dA)(\nabla'\omega-u'\Omega_0\omega)\\
&=& e^{u'dA}\omega-\omega - K((\nabla'-u'\Omega_0)\omega).
\end{eqnarray*}

(b) Let $B_{\frac{1}{2}}=B_0+dA$. It suffices to prove that   the following
diagrams are  commutative:
\begin{equation}
\xymatrix{
  & \qquad \qquad H^*(\Omega^*_c(\Lambda\XX,L')((u)),\nabla'-u'\Omega_0)
     \ar[d]^{Id}\\
HP_*(\aA)\qquad \ar[ur]^{\hspace{-40pt}\tau_{\nabla_0,B_0}}
         \ar[r]_{\hspace{-20pt}\tau_{\nabla_1,B_{\frac{1}{2}}}}
         \ar[dr]_{\tau_{\nabla_1,B_1}}
  &\qquad \qquad H^*(\Omega^*_c(\Lambda\XX,L')((u)),\nabla'-u'\Omega_0)
    \ar[d]^{e^{u'\beta_{}}}\\
  &\qquad \qquad H^*(\Omega^*_c(\Lambda\XX,L')((u)),\nabla'-u'\Omega_1)
}
\end{equation}

For the lower triangle, note that
\begin{eqnarray*}
\lefteqn{e^{u'\beta}*\tilde{a}_0*e^{u's_0 B_{\frac{1}{2}}}*\nabla_1 a_1*\cdots
*\nabla_1a_k* e^{u's_kB_{\frac{1}{2}}}}\\
 &=&
\tilde{a}_0*e^{u's_0\beta}*e^{u's_0 B_{\frac{1}{2}}}*\nabla_1 a_1*\cdots
*\nabla_1a_k *e^{u's_k\beta}*e^{u' s_k B_{\frac{1}{2}}}\\
&=&\tilde{a}_0*e^{u's_0B_{1}}*\nabla_1 a_1*\cdots *
\nabla_1a_k *e^{u's_k B_{1}}, 
\end{eqnarray*}
where we use the relations $s_0+\cdots+s_k=1$ and
$\beta=B_1-B_{\frac{1}{2}}$. Note that since $\del\beta=0$,
 $\beta$  commutes with   every element in $\Omega^*(\gm,L)$
according to Lemma  \ref{lem:star}.
\par\medskip
For the upper triangle, we proceed with the standard argument as in
\cite[Proposition~5.6]{Mat-Ste04} (but the proof is simpler here).
Let $\gm^\sharp=\gm\times [0,1]$, and  $\tilde{\gm}^\sharp
=\tilde{\gm}\times [0,1]$. Then $\tilde{\gm}^\sharp\to \gm^\sharp\toto
M\times [0,1] $ is an $S^1$-central extension. Let
$\theta^\sharp=\theta +t\del A$, which is   considered as
a 1-form on $\tilde{\gm}^\sharp$. It is simple  to see that
$\theta^\sharp$ defines a connection of the $S^1$-central extension
$\tilde{\gm}^\sharp\to \gm^\sharp\toto
M\times [0,1] $. Moreover, since
$$ d\theta^\sharp=d\theta +dt\wedge \del A +t d\del A=
\del (B+dt \wedge A+t dA), $$
and 
$$ d(B+dt \wedge A+t dA)=dB=\Omega, $$
it follows that $B^\sharp=B+dt \wedge A+t dA\in \Omega^{2}(M\times [0, 1])$
is its  curving, and $\Omega$, being considered as an element in
$\Omega^3 (M\times [0, 1])^{\gm^\sharp}$, is the curvature. 
%
Applying  Proposition~\ref{prop:chain-map},
we have
\begin{equation}
\label{eq:tbB}
\tau^\sharp\smalcirc(b+u\B)=(u\nabla^{\sharp'}-2\pi i u^2\Omega)\smalcirc
\tau^\sharp,
\end{equation}
where $\tau^\sharp=\tau_{\nabla^\sharp,B^\sharp}$.
Note that when being restricted to $S \gm^\sharp (\cong S\gm \times [0, 1])$,
$\theta^\sharp$ is equal to $\theta$ since $\del A|_{S\gm }=0$.
Therefore its corresponding covariant differential
  $\nabla^{\sharp'}$ is equal to $\nabla' +(dt) \frac{d}{dt}$, i.e.
for any $v\in T_h S\gm$ and  any $t$-dependent section $s_t\in \gm (L')$
being considered as a section of $L'\times [0, 1]\to S\gm \times [0, 1]$,
we have
$$\nabla^{\sharp'}_{(v, \frac{d}{dt})}s_t =\nabla' s_t+ \frac{ds_t}{dt} .$$
Contracting Eq. \eqref{eq:tbB}
 by the vector field $\frac{d}{dt}$, dividing by $u$ and
integrating over $[ 0 , 1]$, we obtain 
$$\tau_{\nabla_1,B_1}-\tau_{\nabla_0,B_0} = K\smalcirc (b+u\B)-(u
\nabla'-2\pi i u^2\Omega)\smalcirc K, $$
where $K=u^{-1}\int_0^1 \iota_{\frac{d}{dt}}\tau^\sharp\,dt$.
\end{pf}

\subsection{Proof of the main theorem}

Our goal in this  subsection is to prove
Theorem~\ref{thm:chern}. The idea is that $K_\alpha^*(\XX)\otimes\cc$,
$HP_*(C_c^\infty(\gm,L))$ and $H^*(\XX,\alpha)$
agree locally, since an orbifold is locally a crossed-product of
a manifold by a finite group; and each of these cohomology functors
admits Mayer-Vietoris sequences,
hence they agree globally.

For every open subset $U$ of $|\XX|$, we denote by $\gm_U$ the restriction
of $\gm$ to $U$, i.e., 
$$\gm_U =\{\gamma\in \gm |s(\gamma) , t(\gamma)\in \pi^{-1} (U)\}$$
where $\pi: M\to |\XX|$ is the projection. The corresponding
orbifold of $\gm_U$ is denoted by $\UU$.

\begin{lem}\label{lem:local-global}
For every open subset $U$ of $|\XX|$, let
$H^*(U, \alpha)=K_\alpha^*(\UU)\otimes\cc$,
 $HP_*(C_c^\infty(\gm_U,L))$, 
or $H^*_c(\UU,\alpha )((u))=H_c^{*+2\zz}(\UU,\alpha)$ ($*\in\zz/2\zz$).
\begin{itemize}
\item[(a)] If $(U_i)$ is an increasing net of open subsets of $|\XX|$
which covers $|\XX|$, then $H^*(|\XX|,\alpha)=\lim_i H^*(U_i,\alpha)$;
\item[(b)] If $|\XX|$ is covered by two open subsets $U$ and $V$,
then there is a Mayer-Vietoris exact sequence
\begin{equation*}
\xymatrix{
H^0(U\cap V,\alpha)\ar[r]& H^0(U,\alpha)\oplus H^0(V,\alpha)\ar[r]&
 H^0(|\XX|,\alpha)\ar[d]\\
H^1(|\XX|,\alpha)\ar[u] & H^1(U,\alpha)\oplus H^1(V,\alpha)\ar[l]&
 H^1(U\cap V,\alpha)\ar[l]
}
\end{equation*}
\end{itemize}
\end{lem}

\begin{pf}
(a) For $K_*$, this follows from the facts that $K$-theory commutes with
inductive limits, and that $(C^*(\gm_{U_i}, L))$ is an increasing net
of ideals in $C^*(\gm , L)$ whose union is dense in $C^*(\gm , L)$.

For $HP_*$ and $H^*_c(-,\alpha)$, this is obvious.

(b) For $K_*$, this is proven in \cite[Proposition~3.9]{TXL04}.
For $HP_*$, this follows from \cite{CQ97} (see also \cite{Nis97}) and for
$H^*_c(-,\alpha)$, the proof is standard using smooth partitions of unity
(see Lemma~\ref{lem:partition-unity} below).
\end{pf}

\begin{lem}\label{lem:partition-unity}
Given an open cover $\uU=(U_i)$ of $|\XX|$,
there exists a partition of unity subordinate to $\uU$ consisting
of smooth functions on $\XX$.
\end{lem}

\begin{pf}
Let $\pi:M\to |\XX|$ be the  projection map.
For each $i$, choose an open submanifold 
 $V_i$ of $M$ such that $\pi(V_i)=U_i$.
 Let $\gm'\toto \amalg V_i$ be the pull-back of
$\gm\toto M$ via the \'etale  map $\amalg V_i \to M$.
Let $c=(c_i)_{i\in I}$ be a smooth cutoff function for the proper
groupoid $\gm'\toto \amalg V_i$
\cite[Proposition~6.7]{Tu99}. By definition,
$c_i: V_i\to \rr_+$ is a smooth function such that for all $i$ and all
$x\in V_i$,
$\sum_j\sum_{g\in \gm^x} c_j(s(g)) =1$ (by convention,
$c_j$ is extended by zero outside $V_j$).
Let $\varphi_i(x)=\sum_{g\in \gm^x} c_i(s(g))$ ($x\in M$). Then
$\varphi_i$ is clearly $\gm$-invariant, smooth, and $\sum\varphi_i=1$.
\end{pf}

Thus, in order to prove Theorem~\ref{thm:chern},
 it suffices, by induction using a five-lemma argument and  passing
to the inductive
limit, to prove the following

\begin{prop}
\label{pro:local}
For  every $\bar x\in |\XX|$, there exists   an open neighborhood
$U$ such that for every open subset $V\subset U$, the  homomorphisms
\begin{equation}\label{eqn:iso-local}
K_*(C_c^\infty(\gm_{V},L))
\otimes\cc\stackrel{ch}{\to}
HP_*(C_c^\infty(\gm_{V},L))\stackrel{\tau}{\to}
H^*_c(\VV,\alpha)
\end{equation}
are isomorphisms.
\end{prop}

\begin{lem}\label{lem:local}
For each $\bar x\in |\XX|$, there exists an open neighborhood $U$
of $\bar x$ such that $\alpha_{|U} \in H^3(  \XX|_U, \zz )$ is represented by
an $S^1$-central extension:
\begin{equation}\label{eqn:local-central-ext}
S^1\to U'\rtimes \tilde{G}\to U'\rtimes G\toto U',
\end{equation}
 which is the pull back 
  of a group $S^1$-central extension  of finite order: 
\begin{equation}\label{eqn:extension-G}
S^1\to \tilde{G}\to G,
\end{equation}
induced by a $\zz_n$-central extension:
\begin{equation}
\label{eqn:extension-G'}
1\to \exp\left(\frac{2\pi i\zz}{n}\right)\to G'\to G\to 1,
\end{equation}
where $\tilde G = G'\times_{\exp(\frac{2\pi i\zz}{n})} S^1$.
 Here  $U'$ is an Euclidean ball and $G$ is the stabilizer of $\bar x$
acting on $U'$ by isometries.
\end{lem}
\begin{pf}
First take a nice  orbifold chart around $\bar x$ of the form $U'\rtimes G$
as in \cite{Moe-Pro97}, where $G$ is a finite group. 
 Since $U'$ is $G$-contractible, we have $H^3_G(U',\zi)\cong H^3_G(pt,\zi)$.
It thus follows that
 the class $\alpha_{|U}$ is represented by an $S^1$-central extension
of $U'\rtimes G$,
 which is the   pull back from an $S^1$-central extension of
the form~(\ref{eqn:extension-G}). Note that the class
$\alpha'\in H^3(G\lcom,\zi)$ of the extension (\ref{eqn:extension-G})
must be   torsion since the image of $\alpha'$ in $H^3(G\lcom, \rr)=
H^3(G\lcom ,\zi)\otimes\rr$ is zero.

It follows that, if $n$ is such that
$n\alpha'=0$, the extension (\ref{eqn:extension-G})
is obtained from (\ref{eqn:extension-G'}).
\end{pf}

\begin{lem}
\label{lem:torsion}
Assume that  $\gm$ is an   \'etale proper groupoid,
$S^1\to\tilde{\gm}\to\gm\toto M$ is a  torsion class $S^1$-central
extension  induced from a     $\zz_n$-central extension
$\exp{\frac{2\pi i \zz}{n}} \to \tilde\gm'\to\gm\toto M$.
Then there is an induced  flat gerbe connection
on $\tilde{\gm}\to\gm\toto M$.
\end{lem}
\begin{pf}
Let $L=\tilde{\gm}\times_{S^1} \cc$ be the induced line bundle
over $\gm$.  For each $q\geq 0$,  the space  
$\Omega^q (\tgm')$ of $q$-forms on $\gm$ admits a 
decomposition
\begin{equation}
\label{eq:dec}
\Omega^q (\tgm')= \oplus_{k=0}^{n-1} \Omega^q (\gm,L^{\otimes k}),
\end{equation}
where  $\Omega^q (\gm,L^{\otimes k})$ can be  naturally identified with
the space $\{\omega \in \Omega^q (\gm ')|\; z\cdot \omega
=(\exp{\frac{2\pi izk}{n}})\omega, \ \forall z\in \zz_n \}$
 (this is analogue to \cite[Proposition~3.2]{TXL04}), where $z\cdot$
denotes the  induced $\zz_n$-action on $\Omega^q (\tgm')$.
Indeed, if we consider the operator $T$ on $\Omega^q (\tgm')$
given by $T (\omega )= e^{-\frac{2\pi i}{n}}\cdot \omega $,
then Eq. (\ref{eq:dec}) is just the decomposition of
$\Omega^q (\tgm')$ into the  eigenspaces of $T$:
$$\Omega^q (\tgm')=\oplus_{k=0}^{n-1}
\ker (T-e^{\frac{2\pi i k}{n}}{\mathrm{Id}}).$$
It is simple to see that each eigenspace of $T$, i.e.
$\Omega^q (\gm,L^{\otimes k})$ is stable under the de Rham differential
 $d: \Omega^* (\tgm') \to \Omega^{*+1} (\tgm')$. In particular, one
obtains a degree $1$ differential operator $\nabla_1$ 
on $\oplus_q \Omega^q (\gm,L)$, which is easily seen to  satisfy
the axioms of a covariant differential. Therefore there is an
induced  linear connection on $L\to \gm$, which must be flat since
$d^2=0$. To show that $\nabla_1$ satisfies the gerbe connection
condition, we have to check that $\nabla_1(f_1*f_2)=
\nabla_1 f_1 * f_2 + f_1*\nabla_1 f_2$ (see Lemma~\ref{lem:connec-star}),
where $*$ is the convolution product on $C_c^\infty(\gm,L)$.
Let $*'$ be the convolution product on $\tilde{\gm}'$.
Then, identifying $C_c^\infty(\gm,L)$ with a subspace of
$C_c^\infty(\tilde{\gm}')$ as above, we have
$f_1*'f_2=n f_1*f_2$. So it suffices to show that $d(f_1*'f_2)
=df_1 *' f_2 + f_1 *' df_2$. The latter  follows from
Lemma~\ref{lem:connec-star} applied to the groupoid $\tilde{\gm}'$
endowed with the trivial gerbe.
\end{pf}

Note that the above construction works for compactly supported
differential forms as well. In particular,
 $\aA'=C_c^\infty(\tilde\gm',\cc)$ is the direct sum
\begin{equation}\label{eqn:direct-sum}
\aA'=\oplus_{k=0}^{n-1} C_c^\infty(\gm,L^{\otimes k}).
\end{equation}

We thus obtain the following decompositions:
\begin{eqnarray*}
\nonumber
K_*(C_c^\infty(\tilde\gm'))&=&
\oplus_{k=0}^{n-1} K_*(C_c^\infty(\gm,L^{\otimes k}))\\
\nonumber
HP_*(C_c^\infty(\tilde\gm'))&=&
\oplus_{k=0}^{n-1} HP_*(C_c^\infty(\gm,L^{\otimes k})).
\end{eqnarray*}

On the other hand, Lemma \ref{lem:torsion} and its proof imply
that
\begin{equation}
\label{eqn:decomp-omega}
H^{*}( \Omega_c^*(S\tgm')^{\tilde\gm'}((u)))=
\oplus_{k=0}^{n-1} H^{*}
((\Omega^*_c(S\gm,L^{'\otimes k})^\gm((u)), \nabla_1^{'k}),
\end{equation}
where $\nabla_1$ is the connection on $L\to \gm$ as in Lemma \ref{lem:torsion},
and $\nabla_1^k$ its induced connection
on $L^{\otimes k}\to \gm$ while the upscript $'$ stands for
their restrictions to the closed loops $S\gm$.
 We are now ready to prove Proposition
\ref{pro:local}.\\\\\\\\\\
{\bf Proof of Proposition \ref{pro:local}}
Assume that we already have established the following isomorphisms:
\begin{equation}\label{eqn:iso-gm'}
K_*(C_c^\infty(\tilde\gm'))\otimes \cc
\stackrel{ch}{\to} HP_*(C_c^\infty(\tilde\gm'))
\stackrel{\tau}{\to} H^{*}(\Omega^*_c(S\tilde\gm')^{\tilde\gm'}((u))),
\end{equation}
where $\tau$ is the homomorphism constructed
as in Eq. (\ref{eq:tau}) using the trivial connection
 and curving $0$ by considering the groupoid $\tilde{\gm}'$ being
 equipped with the trivial gerbe. More precisely,

\begin{equation}\label{eqn:tau}
\tau(\tilde{a}_0\otimes a_1\otimes\cdots\otimes a_k)
=\frac{1}{k!}Tr(\tilde{a}_0 *da_1*\cdots *da_k).
\end{equation}

 Then the isomorphisms
(\ref{eqn:iso-gm'}) induce isomorphisms on eigenspaces of $T$
(since all maps in (\ref{eqn:iso-gm'}) commute with $T$). Thus
$$K_*(C_c^\infty(\gm,L))\otimes \cc
\stackrel{ch}{\to} HP_* (C_c^\infty(\gm,L))
\stackrel{\tau_{\nabla_1,0}}{\to}
H^*(\Omega^*_c((S \gm, L')^{\gm} ((u)),\nabla'_1 ))$$
are isomorphisms. From
Proposition~\ref{prop:two-connections}, it follows that
$$K_* (C_c^\infty(\gm,L))\otimes \cc \stackrel{ch}{\to} HP_* (C_c^\infty(\gm,L))
\stackrel{\tau_{\nabla,B}}{\to}
H^*(\Omega^*_c(S \gm,L')^{\gm} ((u)),\nabla'_1-2\pi i  u\Omega )$$
are isomorphisms.

It thus remains to prove that (\ref{eqn:iso-gm'}) are indeed isomorphisms.
Therefore we are reduced to the case when $\gm$ is Morita equivalent
to the crossed-product $U\rtimes G$
of a manifold by a finite group and $\alpha=0$. Of course, we may
also assume that  $V=U$, since $\gm_{|V}$ satisfies the same properties.

In this  case, using the two lemmas below,
we can replace $\gm\toto M$ by its Morita equivalent
groupoid $U\rtimes G\toto U$.

\begin{lem}
\label{lem:4.13}
Assume that $\gm\toto M$ is  an  \'etale  proper  groupoid. Let $(U_i)$ be an open
cover of $M$, $M'=\amalg U_i$ and $\gm'=\{(i,g,j)\vert\;
g\in\gm^{U_i}_{U_j}\}$ the pull back of
$\gm$ under the \'etale map $\amalg U_i\to M$.
The following diagram is commutative:

\begin{equation}\label{mor.coh}
\vcenter{\xymatrix{
{HP_*(C_c^\infty(\gm')) }\rrto^{\phi }\dto^{\tau'} &&
 { HP_*(C_c^\infty(\gm))}\dto^{\tau}\\
{H^*(\Omega^*_c(S\gm' )^{\gm'}((u)))}\rrto^\sim &&
 {H^*(\Omega^*_c (S\gm )^\gm((u)))}}}
\end{equation}
%
where $\tau$ and $\tau'$ are defined by Eq. (\ref{eqn:tau}),  and 
$\phi$ is given  by
$$\phi (\tilde{a}_0\otimes \cdots\otimes a_k)
(g_0,\ldots,g_k)=\sum_{i_0,\ldots,i_k}
\tilde{a}_0(i_0,g_0,i_1)\otimes\cdots\otimes
a_n(i_k,g_k,i_0)$$
(note that $g\mapsto a_j(i_j,g,i_{j+1})$
is a smooth compactly supported function on
$\gm^{U_{i_j}}_{U_{i_{j+1}}}$. Therefore it  can be considered as
a smooth compactly supported function on $\gm$).
Moreover, $\phi$ is an isomorphism.
\end{lem}
\begin{pf}
The fact that $\phi$ is an isomorphism, i.e. that $HP_*(C_c^\infty(\gm))$
only depends on the Morita equivalence class of the \'etale proper
groupoid $\gm\toto M$ is standard; e.g. see  \cite{Cra99,Cra-Moe00}.

That the diagram (\ref{mor.coh}) commutes follows from a direct and elementary computation:
for any $h\in S\gm$,
$(\tau\smalcirc \phi )(\tilde{a}_0\otimes \cdots\otimes a_k)(h)$
and $\tau'(\tilde{a}_0\otimes \cdots\otimes a_k)(h)$
are both equal to
$$\frac{1}{k!}\sum_{g\in\gm^{s(h)}}
\sum_{g_0\cdots g_k = g^{-1}hg}
\sum_{i_0,\ldots,i_k}
\tilde{a}_0(i_0,g_0,i_1)da_1(i_1,g_1,i_2)\cdots
da_k(i_k,g_k,i_0).$$
\end{pf}

\begin{lem}
Let $\gm$ and $\gm'$ be  as in Lemma \ref{lem:4.13}. The following
diagram commutes, 
and the  horizontal maps are isomorphisms:
$$\xymatrix{
K_*(C_c^\infty(\gm'))\ar[r]\ar[d]^{ch}
  & K_*(C_c^\infty(\gm))\ar[d]^{ch}\\
HP_*(C_c^\infty(\gm'))\ar[r]
  & HP_*(C_c^\infty(\gm))
}$$
\end{lem}

\begin{pf}
Let $\gm''=\gm\times (I\times I)$, where $I$ is the
index set of $(U_i )$, and  $I\times I$ is the pair groupoid
equipped with the product $(i,j)(j,k)=(i,k)$. The inclusion of $\gm'$
as an open subgroupoid of $\gm''$ induces a $*$-homomorphism
$C_c^\infty(\gm')\to C_c^\infty(\gm'')$. Choosing $i_0\in I$,
the inclusion $g\mapsto (g,i_0,i_0)$ also induces a $*$-homomorphism
$C_c^\infty(\gm)\to C_c^\infty(\gm'')$. Moreover, these $*$-homomorphisms
are Morita equivalences, hence the following diagram
$$\xymatrix{
K_*(C_c^\infty(\gm'))\ar[r]\ar[d]^{ch}
  & K_*(C_c^\infty(\gm''))\ar[d]^{ch}
  & K_*(C_c^\infty(\gm))\ar[d]^{ch}\ar[l]\\
HP_*(C_c^\infty(\gm'))\ar[r]
  & HP_*(C_c^\infty(\gm''))
  & HP_*(C_c^\infty(\gm))\ar[l]
}$$
is commutative and the horizontal arrows are isomorphisms.
\end{pf}

Now, the case $\gm =U\rtimes G$ is covered by Baum and Connes
\cite[Theorem~1.19]{Bau-Con88}. More precisely, the proof of
Proposition \ref{pro:local} and
 Theorem~\ref{thm:chern} will be completed thanks to the following:

\begin{lem}
Let $G$ be a finite group acting on a manifold $M$.
The following homomorphisms
\begin{equation}
\label{eqn:iso-MG}
K_*(C_c^\infty(M\rtimes G))\otimes \cc \stackrel{ch}{\to}
HP_* (C_c^\infty(M\rtimes G)) \stackrel{\tau}{\to} 
H^*(\Omega^*_c (\hat{M})^G((u))),
\end{equation}
are isomorphisms,
where $\hat{M}=\amalg_{g\in G}M^g$, $\Omega_c^*$ denotes the space of
differential forms which are compactly supported in $\hat{M}/G$,
 and $\tau$
is the map defined by Eq. (\ref{eqn:tau}). Moreover, $\tau\smalcirc ch$ is
the Baum-Connes' Chern character \cite{Bau-Con88} and hence is an
isomorphism.
\end{lem}
\begin{pf}
This is a well known   result. However since we cannot locate
it in literature, we sketch a proof below.

Let us treat some special cases first. If $G$ is the trivial group
then $\tau\smalcirc ch$ is the usual Chern character and $\tau$
is Connes-Hochschild-Kostant-Rosenberg's isomorphism
\cite{Con85,Pfl98}.

If $M$ is a point, the maps (\ref{eqn:iso-MG})  become
$$R(G)\otimes\cc\stackrel{ch}{\to} \cc[G]^G\stackrel{\tau}{\to} \cc[G]^G$$
where $\cc[G]^G$ denotes the complex-valued functions on $G$
which are invariant under conjugation. Then $ch$ is the character map;
we will show that $\tau$ is the identity map.
 For this purpose,  let us check that $\tau\smalcirc ch$
 is also the character map.
Let $\pi$ be an irreducible representation of $G$, $\chi$ its
character and $d_\pi$ its dimension. Let $f(g)=\frac{d_\pi\chi(g)}{\# G}$.
Then the corresponding element $P\in \cc[G]=C^*(G)$ is a projection
which corresponds to the element $d_\pi [\pi]\in R(G)$ \cite{Dix77}
(the term $\# G$ comes from the fact that we use the counting measure
instead of the normalized
Haar measure). Moreover, one checks immediately that
$(\tau\smalcirc ch)( d_\pi[\pi])(h)=(\tau\smalcirc  ch([P]))(h)
=(\tau(f))(h)=\sum_g f(g^{-1}hg) =d_\pi \chi(h)$.
\par\medskip

Now the case when $G$ acts trivially on $M$ follows from the isomorphisms
$K_*(C_c^\infty(M)\rtimes G)=K_*(C_c^\infty(M))\otimes R(G)$
and the analogue isomorphisms for $HP_*$ and  $H^*_c$.
\par\medskip

Let us turn to the general case.
Let $\sigma(M,G)=\max \{\#{\mathrm{stab}}(x)|\; x\in M\}$.
We proceed by induction on $\sigma(M,G)$. If $\sigma(M,G)=1$, then
$G$ acts freely  on $M$. After replacing $M\rtimes G$
by the Morita equivalent
manifold $M/G$ we are reduced to Connes-Hochschild-Kostant-Rosenberg's 
theorem
as above.

Suppose the proposition is proven for $\sigma(M,G)<N$ and
assume $\sigma(M,G)=N$. Let $U=\{x\in M|\;\#{\mathrm{stab}}(x)<N\}$
and $F=M-U$. Then $U$ is an open invariant subset of $M$ such that
$\sigma(U,G)<N$ and $F$ is a closed submanifold of $M$.
By the induction assumption, the proposition is true for $U$,
so by using
six-term exact sequences associated to the pair $(U,F)$,
we just have to show the proposition for the space $F$. I.e.
we are reduced to the case when all stabilizers have cardinality $N$.
Now, let $M^H=\{x\in M|\; {\mathrm{stab}}(x)=H\}$ and
let $M_H$ be the saturation of $M^H$. Choose a representative $H_i$
for each conjugacy class of subgroups of cardinal $N$.
Since the stabilizer of a point in $M$ is conjugate to one
of the $H_i$'s, $M=\amalg M_{H_i}$ is a
finite partition of $M$. Moreover, since all stabilizers have
cardinality $N$, we have
$M^H=\{x\in M\vert\; {\mathrm{stab}}(x)\supseteq H\}$.
Hence $M^H$ is closed for all $H$. It follows that 
$M=\amalg M_{H_i}$ is a
partition of $M$ into finitely many open and closed subsets,
so we can assume that $M=M_H$ for some $H$.
Then since $M^H\cdot g= M^{g^{-1}Hg}$, we see that
$M\rtimes G$ is Morita equivalent to
$M^H\rtimes K$, where $K$ is the normalizer of $H$. Hence we can assume that,
after replacing $M$ by $M^H$ and $G$ by $K$,
 $H$ is normal and  ${\mathrm{stab}}(x)=H$ for all $x$,
i.e. the action of $G$ comes from a free action of $G/H$.
\par\medskip

Now, for all $x\in M/G=M/(G/H)$, there exists a neighborhood
$V$ of $x$ such that $\pi^{-1}(V)$ is equivariantly diffeomorphic
to $V\times (G/H)$. Cover $M/G$ by such open subsets $V_i$.
Using Mayer-Vietoris exact sequences and an induction argument,
we can assume that $M=V\times (G/H)$. Therefore $M\rtimes G$
is Morita equivalent to the crossed-product $V\times H$ of a manifold
by a trivial group action, which is the case  considered earlier.
\end{pf}


\small
\begin{thebibliography}{00}
\bibitem{Ade-Rua03}
Adem, A., and  Ruan, Y., Twisted orbifold $K$-theory,
\emph{Commun. Math. Phys.} {\bf 237}  (2003),   533--556.

\bibitem{BBM}
 Baum, P., Brylinski, J.-L., and MacPherson, R.,
 Cohomologie \'equivariante d\'elocalis\'ee
{\em  C. R. Acad. Sci. Paris Serie I Math.} {\bf  300} (1985),   605--608.


\bibitem{Bau-Con88}
Baum, P.,  and Connes, A.,
Chern character for discrete groups,
\emph{A f\^ete of topology}, 163--232,
\emph{Academic Press, Boston, MA}, 1988.

\bibitem{BX1}
Behrend, K., and Xu, P.,
$S^1$-bundles and gerbes over differentiable stacks,
\emph{C. R. Acad. Sci. Paris, S\'erie I}  {\bf 336}
(2003), 163-168.


\bibitem{BX2}
Behrend, K.,  and Xu, P., Differentiable stacks and gerbes,  in preparation.

\bibitem{BG94}
Block, J.,  and Getzler, E.,
Equivariant cyclic homology and equivariant differential forms,
\emph{Ann. Sci. ENS} 4e s\'erie, {\bf 27} (1994), 493--527.

\bibitem{Brylinski}
Brylinski, J.-L.,
Loop spaces, characteristic classes and geometric quantization,
{\em Progress in Mathematics} {\bf 107},
 Birkh\"auser, Boston, MA, 1993.

 
            
 

\bibitem{Bry-Nis94}
Brylinski, J.-L., and  Nistor, V.,
Cyclic cohomology of \'etale groupoids,
\emph{$K$-Theory} {\bf 8}  (1994),  341--365. 

\bibitem{Bur85}
Burghelea, D.,
The cyclic homology of the group rings,
\emph{Comment. Math. Helv.} {\bf 60}  (1985),   354--365.

\bibitem{Con85}
Connes, A., {Non-commutative Differential Geometry},
\emph{Publ. Math. IHES} {\bf 62} (1985) 257--360.

\bibitem{Cra99}
Crainic, M.,  Cyclic cohomology of \'etale groupoids:
the general case,  \emph{$K$-Theory} {\bf 17}  (1999),   319--362.

\bibitem{Cra03}
Crainic, M.,
Differentiable and algebroid cohomology, van Est isomorphisms, and
characteristic classes, 
{\em Comment. Math. Helv.} {\bf 78} (2003),  681--721.

\bibitem{Cra-Moe00}
Crainic, M., and  Moerdijk, I. A homology theory for \'etale groupoids,
\emph{J. Reine Angew. Math.} {\bf 521}  (2000), 25--46.


\bibitem{CQ97}
Cuntz, J.,  and Quillen, D.,
Excision in bivariant periodic cyclic cohomology,
\emph{Invent. Math.} {\bf 127} (1997), 67-98.

\bibitem{Dix77}
Dixmier, J.,
$C^*$-algebras,
North Holland Publishing Co.,
Amsterdam--{New York}--Oxford, 1977.


\bibitem{Fei-Tsy87}
Feigin, B., and  Tsygan, B.,
Additive $K$-theory,
$K$-theory, arithmetic and geometry (Moscow, 1984--1986),  67--209,
\emph{Lecture Notes in Math.}, {\bf 1289}, Springer, Berlin, 1987.

\bibitem{giraud}
Giraud, J.,
Cohomologie non ab\'elienne,
{\em Grundlehren der mathematischen Wissenschaften} {\bf 179}
{Springer-Verlag}, Berlin, 1971.


\bibitem{Gor99}
 Gorokhovsky, A.,
{Characters of cycles, equivariant characteristic
classes and Fredholm modules},
\emph{Commun. Math. Phys.} {\bf 208} (1999) 1--23.

\bibitem{KMRW98}
Kumjian, A.,  Muhly, P.,  Renault, J. and Williams, D.,  
The Brauer group of a locally compact
groupoid, \emph{Amer. J. Math.} {\bf  120} (1998),  901--954.

\bibitem{LTX}
Laurent-Gengoux, C., Tu, J.-L.  and Xu, P.,
Chern-Weil map for principal bundles over groupoids, 
{\em Math. Z.} (to appear), {\tt math.DG/0401420}.

\bibitem{Lod98}
Loday, J. L.
Cyclic homology,
{\em Grundlehren der Mathematischen Wissenschaften}
{\bf 301} Springer-Verlag, Berlin, 1998.

\bibitem{LU02}
Lupercio, E.,  and Uribe, B.,
Loop groupoids, gerbes, and twisted sectors on orbifolds,
{Orbifolds in mathematics and physics} (Madison, WI, 2001), 163--184,
{\em Contemp. Math.} {\bf 310}, Amer. Math. Soc., Providence, RI, 2002.



\bibitem{Mat-Ste04}
Mathai, V., and Stevenson, D.,
On a generalized Connes-Hochschild-Kostant-Rosenberg theorem,
 \emph{Adv. Math} (to appear). {\tt math.KT/0404329}

\bibitem{Moe02}
Moerdijk, I.,  Orbifolds as groupoids: an introduction. 
Orbifolds in mathematics and physics (Madison, WI, 2001),  205--222,
{\em Contemp. Math.} {\bf 310}, Amer. Math. Soc., Providence, RI, 2002.

\bibitem{Moe-Pro97}
Moerdijk, I. and  Pronk, D. A., Orbifolds, sheaves and groupoids,
{\em $K$-Theory}  {\bf 12}  (1997),  3--21.

\bibitem{Nis90}
Nistor, V.,
Group cohomology and the cyclic cohomology of crossed products,
\emph{Invent. Math.} {\bf  99}  (1990),  411--424.

\bibitem{Nis97}
Nistor, V.,
Higher index theorems and the boundary map in cyclic cohomology,
\emph{Doc. Math} {\bf 2} (1997) 263--295.

\bibitem{Pfl98}
Pflaum, M.,
 { On continuous Hochschild homology and cohomology
groups}, \emph{Lett. Math. Phys.} {\bf 44} (1998), 43--51.

\bibitem{Quillen}
 Quillen, D.,
 Algebra cochains and  cyclic cohomology,
{\em Publ. Math. IHES} {\bf 68} (1989), 139-174.



\bibitem{Rua03}
Ruan, Y.,
Discrete torsion and twisted orbifold cohomology,
\emph{J. Sympl. Geom.} {\bf  2} (2003), 1--24.

\bibitem{Tu99}
Tu, J.-L.,  La conjecture de Novikov pour les feuilletages hyperboliques,
\emph{$K$-Theory} {\bf 16}  (1999),   129--184.

\bibitem{TXL04}
Tu, J.-L., Xu, P.,  and Laurent-Gengoux, C.,
Twisted $K$-theory of differentiable stacks,
{\em Ann. Sci. ENS.} {\bf 37} (2004), 841-910.

\bibitem{Witten}
Witten, E.,
$D$-branes and $K$-theory,
\emph{J. High Energy Phys.} {\bf 12} (1998),
19-44.

\end{thebibliography}
\end{document}